\documentclass{article}

\usepackage{arxiv}
\usepackage{graphicx}
\usepackage{subfigure}
\usepackage{stmaryrd}
\usepackage{algorithm}
\usepackage{algorithmic}
\usepackage{xcolor}
\usepackage{bbm}
\usepackage{enumitem}
\usepackage{times,amsmath}
\usepackage{latexsym,amssymb,lscape,multirow}
\usepackage{bm}
\usepackage[numbers]{natbib}
\allowdisplaybreaks[4]

\usepackage{tikz}
\usepackage{pgfplots}
\usetikzlibrary{shapes,decorations}
\usetikzlibrary{fit}

\usepackage{dcolumn,amsthm}

\renewenvironment{proof}[1][Proof]{\noindent\textit{#1. } }{\hfill$\square$}

\newcommand{\R}{\mathbb{R}}
\newcommand{\N}{\mathbb{N}}
\newcommand{\Rn}{\mathbb{R}^n}
\DeclareMathOperator*{\argmin}{argmin}

\numberwithin{equation}{section}

\begin{document}

\title{A Trust-Region Method for Nonsmooth Nonconvex Optimization}
% Short title for running heads:
%\shorttitle{A Trust-Region Method For Nonsmooth Nonconvex Optimization}

\author{%
{\sc
Ziang Chen\thanks{Email: ziang@math.duke.edu}} \\[2pt]
Department of Mathematics, Duke University, USA\\[6pt]
{\sc and}\\[6pt]
{\sc Andre Milzarek\thanks{Email: andremilzarek@cuhk.edu.cn. A$.$ Milzarek is partly supported by the Fundamental Research Fund -- Shenzhen Research Institute for Big Data (SRIBD) Startup Fund JCYJ-AM20190601.}}\\[2pt]
School of Data Science (SDS), The Chinese University of Hong Kong, Shenzhen, \\ Shenzhen Research Institute for Big Data (SRIBD), CHINA\\[6pt]
{\sc and}\\[6pt]
{\sc
Zaiwen wen\thanks{Email: wenzw@pku.edu.cn. Z. Wen is partly supported by the NSFC grant 11831002, and the Beijing Academy  of Artificial Intelligence.}} \\[2pt]
Beijing International Center for Mathematical Research, Peking University, CHINA}

% Short list of authors for running heads:
%\shortauthorlist{Z. Chen, A. Milzarek, and Z. Wen}

\maketitle

\begin{abstract}
% Body of abstract:
{We propose a trust-region type method for a class of nonsmooth nonconvex optimization problems where the objective function is a summation of a (probably nonconvex) smooth function and a (probably nonsmooth) convex function. The model function of our trust-region subproblem is always quadratic and the linear term of the model is generated using abstract descent directions. Therefore, the trust-region subproblems can be easily constructed as well as efficiently solved by cheap and standard methods. When the accuracy of the model function at the solution of the subproblem is not sufficient, we add a safeguard on the stepsizes for improving the accuracy. For a class of functions that can be ``truncated'', an additional truncation step is defined and a stepsize modification strategy is designed. The overall scheme converges globally and we establish fast local convergence under suitable assumptions. 
In particular, using a connection with a smooth Riemannian trust-region method, we prove local quadratic convergence for partly smooth functions under a strict complementary condition. Preliminary numerical results on a family of $\ell_1$-optimization problems are reported and demonstrate the efficiency of our approach.}

% Keywords:
\end{abstract}

 \section{Introduction}
We consider the unconstrained nonsmooth nonconvex optimization problems of the composite form:
 \begin{equation} \label{obj fun:composite}
 \min_{x\in \mathbb{R}^n}\psi(x):=f(x)+\varphi(x),
 \end{equation}
 where $f:\mathbb{R}^n\rightarrow \mathbb{R}$ is a continuously
 differentiable but probably nonconvex function and
 $\varphi:\mathbb{R}^n\rightarrow\mathbb{R}$ is real-valued and convex. The composite program \eqref{obj fun:composite} is a special form of the general nonsmooth nonconvex optimization problems
 \begin{equation}
     \min_{x\in \mathbb{R}^n}\psi(x),\label{obj fun}
 \end{equation}
 where the objective function $\psi : \Rn \to \R$ is locally Lipschitz continuous, and has numerous applications, such as $\ell_1$-regularized problems \citep{Tibshirani-96, Shevade-03,Don06,KohKimBoy07,BecBobCan11}, group sparse problems \citep{CotRaoEngKre05,YuaLin06,Meier-08,SunLiuCheYe09}, penalty approaches \citep{BacJenMaiObo11}, dictionary learning \citep{mairal2009online,Gangeh-15}, and matrix completion \citep{CanRec09,CaiCanShe10,Hastie-15}.

 \subsection{Related Work} Different types of nonsmooth trust-region methods have already been proposed and analyzed for the general optimization problem \eqref{obj fun} throughout the last two decades. Several of these nonsmooth trust-region methods utilize abstract model functions on a theoretical level which means that the model function is typically not specified. In \citep{Nonsmooth-TR-Dennis-95}, a nonsmooth trust-region method is proposed for \eqref{obj fun} under the assumption that $\psi$ is regular. A nonsmooth trust-region algorithm for general problems is investigated in \citep{Nonsmooth-TR-Qi-Sun-94}. In this work, an abstract first-order model is considered that is not necessarily based on subgradient information or directional derivatives. Extending the results in \citep{Nonsmooth-TR-Qi-Sun-94}, Grohs and Hosseini propose a Riemannian trust-region algorithm, see \citep{Nonsmooth-TR-Riemannian-Grohs-16}. Here, the objective function is defined on a complete Riemannian manifold. All mentioned methods derive global convergence under an assumption similar to the concept of a ``strict model'' stated in \citep{Cutting-Plane-Oracles-Noll-10}. Using this concept, a nonsmooth bundle trust-region algorithm with global convergence is constructed in \citep{Nonsmooth-TR-Apkarian-16}.
 
In \citep{Nonsmooth-TR-Christof-18}, a hybrid approach is presented using simpler and more tractable 
 quadratic model functions. The method switches to a complicated second model if
 the quadratic model is not accurate enough and if it is strictly necessary. In \citep{TR-quadratic-model-Akbari-15},
 a quadratic model function is analyzed where the first-order term is derived from a suitable approximation of the steepest descent direction and the second-order term is updated utilizing a BFGS scheme. The authors apply an
 algorithmic approach proposed in \citep{appro-steepest-descent-Amiri-12} to compute the
 approximation of the $\epsilon$-subdifferential and steepest descent direction. Another class of methods employs smoothing techniques. In \citep{TR-no-derivative-Garmanjani-16}, the authors first present a smooth trust-region method without using derivatives, and then, in the nonsmooth case, use this methodology after smoothing the objective function. Furthermore, trust-region algorithms for nonsmooth
 problems can be developed based on smooth merit functions for the problem. In \citep{Nonsmooth-convex-TR-Sagara-05}, a nonsmooth convex
 optimization is investigated and the Moreau envelope is considered as a smooth merit
 function. A smooth trust-region method is performed on the smooth merit function, where the second-order term of the model function is again updated by the BFGS formula.
 
 Bundle methods are an important and related class of methods for nonsmooth problems \citep{Bundle-Lemarechal-95, Survey-Bundle-Makela-92, Survey-Bundle-Makela-02, Compare-nonsmooth-methods-Karmitsa-12, Bundle-nonsmooth-Lemarechal-94, Bundle-Lan-15, TR-Bundle-Kiwiel-89}. The ubiquitous cutting-plane model in bundle methods is polyhedral, i.e., the supremum of a finite affine family. This model builds approximations of convex functions based on the subgradient inequality. In \citep{bundle-schramm-92}, an efficient bundle technique for convex optimization has been proposed; in \citep{bundle-Oliveira-14}, a convex bundle method is derived to deal with additional noise, i.e., the case when the objective function and the subgradient can not be evaluated exactly. Different modifications of the bundle ideas for nonconvex problems have been established in \citep{Cutting-Plane-Oracles-Noll-10, bundle-schramm-92}. In \citep{bundle-Hare-16} and \citep{bundle-Noll-13}, the authors consider bundle methods for nonsmooth nonconvex optimization when the function values and the subgradients of $\psi$ can only be evaluated inexactly.
 
Local convergence properties and rates for nonsmooth problems are typically studied utilizing additional and more subtle structures. In this regard, some fundamental and helpful concepts are the idea of an ``active manifold'' and the family of ``partly smooth'' functions introduced by \citep{Lewis-02}. In particular, the problem \eqref{obj fun:composite} has been investigated when the nonsmooth term $\varphi$ is partly smooth relative to a smooth active manifold. The so-called finite activity identification is established for forward-backward splitting methods by \citep{partly-smooth-FBE-Liang-17} and, more recently, for SAGA/Prox-SVRG by \citep{partly-smooth-SAGA-Poon-18}. After the identification, those algorithms enter a locally linear convergence regime. In \citep{Finite-active-identify-Hare-04}, the authors use partial smoothness and prox-regularity to identify the active constraints after finitely many iterations, which is an extension of other works on finite constraint identification, see \citep{Finite-active-identify-Burke-88, Finite-active-identify-Burke-90, Finite-active-identify-Wright-93}. After identifying the active manifold, the nonsmooth problem may become a smooth optimization on a Riemannian manifold. Some algorithms and analysis for Riemannian trust-region methods were studied in the literature; see e.g., \cite{absil-08, Absil-07, Huang_RieTR, Baker_RieTR, Baker_RieTR_thesis}. There are numerous applications of Riemannian trust-region methods, such as eigenproblems \cite{Baker_RieTR_thesis, Baker_RieTR_eigen}, low-rank matrix completion \cite{Boumal_RieTR}, and tensor problems \cite{Breiding_RieTR}.
 
We note that for composite programs, nonsmooth trust-region methods have also been studied in the literature, such as \cite{Fletcher_TR_Composite, Fletcher_TR_Composite2, Yuan_TR_Composite, Yamakawa_TR_Composite, Grapiglia_TR_Composite} for $\min\ g(x)+h(f(x))$ where $f$ and $g$ are smooth and $h$ is convex, \cite{Yuan_TR_Composite2, Bannert_TR_Composite} for $\min\ h(f(x))$ where $f$ is smooth and $h$ is convex, \cite{De_TR_Composite} for $\min\ h(f(x))$ where $f$ is locally Lipschitz and $h$ is smooth and convex and \cite{Burke_TR_Composite} for $\min\ g(x)+h(f(x))$ where $f$ is smooth and $g$ and $h$ are convex. For the problem \eqref{obj fun:composite}, there are also other efficient methods, such as gradient-type methods \citep{FukMin81,Nesterov-07}, semismooth Newton methods, \citep{Semismooth-L1-Milzarek-14, Semismooth-Lasso-Li-18, Semismooth-composite-Xiao-18}, proximal Newton methods \citep{Patrinos-13, Prox-Newton-Lee-14, Sto-prox-Newton-Wang-17}, or forward-backward envelope-based (quasi-)Newton methods \citep{FB-Newton-Patrinos-14,FB-Newton-Stella-17}. %, and stochastic model-based methods \citep{DuchiJohn-18}.

More developments about trust-region methods may be found in review papers such as \cite{Yuan_review}.

 \subsection{Our Contribution} In this work, we propose and investigate a trust-region method for nonsmooth composite programs. The approach utilizes quadratic model functions to approximate the underlying nonsmooth objective function. This methodology leads to classical and tractable trust-region subproblems that can be solved efficiently by standard optimization methods if the second-order information is symmetric. We also discuss an efficient subproblem solver in the case that the second-order information does not stem from a symmetric matrix. The linear part of our proposed quadratic model can be based on the steepest descent direction or other directions such as proximal gradient-type descent directions. Our algorithm contains the following steps: computation of a model function and (approximate) solution of the associated trust-region subproblem; acceptance test of the calculated step; determination of a suitable stepsize by a cheap method followed by some stepsize safeguards and a second acceptance test (if the first test is not successful); update of the trust-region radius and a modification step via a novel truncation mechanism.

 In order to control the approximation error between the quadratic model function and the nonsmooth objective function, we define a stepsize safeguard strategy that tries to avoid points along a specific direction at which the directional derivative is not continuous. Specifically, this strategy tries to guarantee that the objective function is directionally differentiable along a specific direction. Since a direct implementation of such a strategy can yield arbitrarily small stepsizes, we consider functions which can be truncated and propose an additional truncation step that allows to enlarge the stepsize. This modification is an essential and new part in our global convergence theory. We verify that the family of functions that can be truncated is rich and contains many important examples, such as the $\ell_1$-norm, $\ell_{\infty}$-norm or group sparse-type penalty terms. Moreover, we provide a detailed global convergence analysis of the proposed trust-region framework. In particular, we show that every accumulation point of a sequence generated by our algorithm is a stationary point. Global convergence of nonsmooth trust-region methods typically requires a certain uniform accuracy assumption on the model which coincides with the concept of the already mentioned strict model proposed by \citep{Cutting-Plane-Oracles-Noll-10}. Our assumptions are similar to these standard requirements and can be verified for a large family containing polyhedral problems and group lasso. Furthermore, we also show how a strict model -- aside from utilizing the original objective function -- can be constructed.
 
 We analyze the local properties of the nonsmooth trust-region method for \eqref{obj fun:composite} when $\varphi$ is a partly smooth function. In particular, it is possible to establish quadratic convergence of our approach in this case. We assume that the underlying manifold is an affine subspace and that a strict complementary condition holds. After the finite activity identification, we transfer our problem to a smooth problem in the affine subspace by proving that an appropriate choice of the first-order and the second-order model coincides with the Riemannian framework. Results from Riemannian trust-region theory can then be applied to derive local quadratic convergence. Additionally, if the nonsmooth term is polyhedral, it can be shown that the Riemannian Hessian can be computed without knowing the underlying manifold.

 \subsection{Organization} The rest of this paper is organized as follows. In Section \ref{sec: prep}, descent directions and several properties of $\psi$ and $\varphi$ used in our algorithm are discussed. In Section \ref{sec: alg}, we present the nonsmooth trust-region framework. In Section \ref{sec: glo_conv}, the global convergence of our method is established. In Section \ref{sec: loc_conv}, we show fast local convergence by studying the nonsmooth composite program for partly smooth $\varphi$. Some preliminary numerical experiments are presented in Section \ref{numerics}.

 \section{Descent Directions and Truncation Operators}
 \label{sec: prep}
In classical trust-region methods, global convergence is established under fairly mild assumptions on the second-order term of the model while the first-order model typically needs to capture the whole gradient information of the objective function, see, e.g., \citep{Num-Opti-Nocedal-06} and the references therein. This underlines the importance of the first-order information in the trust-region method. In this section, we analyze properties of $\psi$ and $\varphi$ as preparation for the construction of suitable linear first-order models.
 
 \subsection{Preliminaries}
 In this work, the expression $\Vert \cdot\Vert$ denotes the $\ell_2$-norm and
 the Frobenius norm for vectors and matrices, respectively. For
 $x\in\mathbb{R}^n$ and $r>0$, $B_r(x):=\{y\in\mathbb{R}^n:\Vert y-x\Vert<r\}$
 denotes the open ball with radius $r$ around $x$. Let $\Lambda$ be a given symmetric positive definite matrix. The proximal operator is defined via
 \begin{equation*}
     \mathrm{prox}_{\varphi}^{\Lambda}(z)= \mathop{\argmin}_{y\in\mathbb{R}^n}\varphi(y)+\frac{1}{2}\Vert y-z\Vert_{\Lambda}^2,
 \end{equation*}
 where $\Vert x\Vert_{\Lambda}^2:=x^T \Lambda x$. We also slightly abuse the notation and write $\mathrm{prox}_{\varphi}^{\lambda}=\mathrm{prox}_{\varphi}^{\Lambda}$ for $\Lambda=\lambda I$. 
 
 The \emph{directional derivative} of a function $h:\mathbb{R}^n\rightarrow\mathbb{R}$ at $x$ along $d$ is denoted by 
$$h'(x;d):=\lim_{t\rightarrow 0^+}\frac{h(x+t d)-h(x)}{t},$$
if it exists. In the composite case $\psi=f+\varphi$ with smooth $f$ and real-valued convex $\varphi$, the directional derivative $\psi'(x;d)$ is well-defined for all $x,d\in\mathbb{R}^n$ and the subdifferential of $\psi$ is defined via
\begin{equation*}
    \partial \psi(x)=\{\nabla f(x)\}+\partial\varphi(x)=\{a\in\mathbb{R}^n\mid \langle a,d\rangle\leq \psi'(x;d),\ \forall\  d\in\mathbb{R}^n\},
\end{equation*}
where $\partial\varphi(x)$ is the usual subdifferential of a convex function. The \emph{steepest descent direction} of $\psi$ is defined as $$d_s(x) :=\begin{cases}\mathop{\argmin}_{\Vert d\Vert\leq 1}\psi'(x;d) & \text{if }0\notin\partial \psi(x),\\
 0 & \text{if } 0\in\partial \psi(x). \end{cases}$$ In this paper, we will repeatedly work with the following normalization condition for a direction $d(x)$:
 
\begin{equation}
 \Vert d(x)\Vert=\begin{cases}1 & \text{if } 0\notin\partial \psi(x),\\
 0 & \text{if } 0\in\partial \psi(x).\end{cases}\label{d(x)}
 \end{equation}
 
We can see that $d_s(x)$ satisfies the property \eqref{d(x)}. We say that a point $x^*$ is a \emph{stationary point} of problem \eqref{obj fun} if $\psi'(x^*;d_s(x^*))=0$, i.e., if and only if, $0\in\partial \psi(x^*)$.

\subsection{Descent Directions}
\label{sec: descent directions}
If the objective function $\psi$ is smooth, then the first-order information is carried in the gradient $\nabla \psi$. We notice that the gradient is directly connected to the steepest descent direction $d_s(x) = -\nabla \psi(x) /\|\nabla \psi(x)\|$ and thus, we can use the following representation $-\nabla \psi(x)=\psi'(x;d_s(x))d_s(x)$. Motivated by this observation, a natural choice of the first-order information and extension in the nonsmooth case is $g(x)=\psi'(x;d_s(x))d_s(x)$.
Since in some cases it might be hard or expensive to directly compute $g(x)=\psi'(x;d_s(x))d_s(x)$, we can also utilize a general descent direction $d(x)$ satisfying \eqref{d(x)} instead. In our model function, we will work with directions of the form $g(x)=u(x)d(x)$ where $d(x)$ is a descent direction, $\psi'(x;d(x)) < 0$, satisfying \eqref{d(x)} and $u(x)$ is an upper bound of $\psi'(x,d(x))$ with
 \begin{equation}\label{u(x)}
 u(x)\begin{cases}
 \in[\psi'(x,d(x)),0)& \text{if } 0\notin\partial \psi(x),\\
 =0 & \text{if } 0\in\partial \psi(x).
 \end{cases}
 \end{equation}
 
This implies $g(x)=0$ if and only if $0\in\partial \psi(x)$, \textit{i.e.}, if $x$ is a stationary point of \eqref{obj fun}. The direction $g(x)$ plays a similar  role as the gradient in the smooth case. We would call $g(x)$ as pseudo-gradient. Our aim in the rest of this subsection is to propose several strategies in the settings of composite programs \eqref{obj fun:composite} for computing and choosing the functions $u(x)$ and $d(x)$.

\subsubsection{Steepest Descent Direction} 
We first compute and express 
\begin{equation}
    g(x)=\psi'(x;d_s(x))d_s(x),\label{g(x):1}
\end{equation}
via the so-called \emph{normal map} \citep{Robinson-92}:
 \begin{equation}
· F_{\mathrm{nor}}^{\Lambda}(z):=\nabla f(\mathrm{prox}_{\varphi}^{\Lambda}(z))+\Lambda(z-\mathrm{prox}_{\varphi}^{\Lambda}(z)),\label{def: normal map}
 \end{equation}
where $\Lambda$ denotes a symmetric and positive semidefinite matrix. We also use the notation $F_{\mathrm{nor}}^{\lambda}=F_{\mathrm{nor}}^{\Lambda}$ in the case $\Lambda=\lambda I$. The next lemma establishes a relation between $d_s(x)$, $F_{\mathrm{nor}}^{\Lambda}(z)$, and $\partial\psi(x)$.

\newtheorem{lem: partial psi}{Lemma}[section]
\begin{lem: partial psi}\label{lem: partial psi}
Let $x\in\mathbb{R}^n$ be given. It holds that

\begin{itemize}
    \item[{(i)}] The direction $d_s(x)$ and the derivative $\psi'(x;d_s(x))$ can be represented as follows: $\psi'(x;d_s(x))=-\mathrm{dist}(0,\partial\psi(x)):=-\min_{v\in\partial\psi(x)}\Vert v\Vert$ and
    \begin{equation*}
    \begin{split}
    d_s(x)=\begin{cases}
        -\frac{\mathbf{P}_{\partial\psi(x)}(0)}{\Vert \mathbf{P}_{\partial\psi(x)}(0)\Vert} & \text{if } 0\notin\partial\psi(x),\\
        0 & \text{if } 0\in\partial\psi(x),
        \end{cases}
    \end{split}  
    \end{equation*}
    where $\mathbf{P}_{\partial\psi(x)}$ denotes the orthogonal projection onto the convex, closed set $\partial\psi(x)$.
    
    \item[{(ii)}] We have $\psi'(x;d_s(x))d_s(x)\in\partial\psi(x)$.

    \item[{(iii)}] $\partial\psi(x)=\{F_{\mathrm{nor}}^{\Lambda}(z): \mathrm{prox}_{\varphi}^{\Lambda}(z)=x\}$.
\end{itemize}
\end{lem: partial psi}

\begin{proof}
(i) Using Fenchel-Rockafellar duality, see \cite[Theorem 15.23]{convex-analysis-hilbert-Bauschke-11}, and the conjugation result $(\iota_{B_{\Vert\cdot\Vert}(0,1)})^*(d)=\sigma_{B_{\Vert\cdot\Vert}(0,1)}(d)=\Vert d\Vert$, we obtain
\begin{equation*}
\begin{split}
\psi'(x;d_s(x))&=\min_d \sigma_{\partial\psi(x)}(d)+\iota_{B_{\Vert\cdot\Vert}(0,1)}(d)\\
    &=-\min_v\iota_{\partial\psi(x)}(v)+\Vert v\Vert=-\text{dist}(0,\partial\psi(x)).
\end{split}
\end{equation*}
The unique solution of the dual problem is given by $v=\mathbf{P}_{\partial\psi(x)}(0)$. By \cite[Corollary 19.2]{convex-analysis-hilbert-Bauschke-11}, the set of primal solutions can be characterized via $d_s(x)\in N_{\partial\psi(x)}(v)\cap\partial\Vert\cdot\Vert(-v)$. Here, the set $N_{\partial\psi(x)}(v):=\{h:\langle h,y-v\rangle\leq 0,\ \forall\  y\in \partial\psi(x)\}$ is the associated normal cone of $\partial\psi(x)$ at $v$. In the case $0\notin\partial\psi(x)$, we have $\Vert v\Vert\neq 0$ and hence, $\partial\Vert\cdot\Vert(-v)=\{-v/\Vert v\Vert\}$. Moreover, since $v$ is a solution of the problem $\min_{y\in\partial\psi(x)}\frac{1}{2}\Vert y\Vert^2$, it satisfies the optimality condition $\langle v,y-v\rangle\geq0,\ \forall\  y\in\partial\psi(x)$. This implies $\{-v/\Vert v\Vert\}\in N_{\partial\psi(x)}(v)$ and $d_s(x)=- \mathbf{P}_{\partial\psi(x)}(0) / {\Vert \mathbf{P}_{\partial\psi(x)}(0)\Vert}$. In the case $0\in\partial\psi(x)$, we have $d_s(x)=0$ by definition.

(ii) Noticing $\psi'(x;d_s(x))=-\Vert \mathbf{P}_{\partial\psi(x)}(0)\Vert$, it follows $\psi'(x;d_s(x))d_s(x)=\mathbf{P}_{\partial\psi(x)}(0)\in\partial\psi(x).$

(iii) By the definition of $\mathrm{prox}_{\varphi}^{\Lambda}(z)$, it can be shown that
\begin{equation}
x=\mathrm{prox}_{\varphi}^{\Lambda}(z) \quad \Longleftrightarrow \quad \Lambda(z-x)\in\partial\varphi(x).\label{prox opti con}
\end{equation}
If $\mathrm{prox}_{\varphi}^{\Lambda}(z)=x$, by \eqref{prox opti con}, we have $F_{\mathrm{nor}}^{\Lambda}(z)=\nabla f(x)+\Lambda(z-x)\in\partial \psi(x)$. If $v\in \partial\psi(x)$, set $z=x+\Lambda^{-1}(v-\nabla f(x))$, then we have $\Lambda(z-x)=v-\nabla f(x)\in\partial\varphi(x)$. According to \eqref{prox opti con}, it holds $x=\mathrm{prox}_{\varphi}^{\Lambda}(z)$. Thus, we obtain $v=F_{\mathrm{nor}}^{\Lambda}(z)$.
\end{proof}

From Lemma~\ref{lem: partial psi} we can immediately derive the following corollary which uses $F_{\mathrm{nor}}^{\Lambda}(z)$ to describe the first-order optimality conditions.

\newtheorem{opti con noraml map}[lem: partial psi]{Corollary}
\begin{opti con noraml map}
A point $x^*\in\mathbb{R}^n$ is a stationary point of problem \eqref{obj fun:composite}, if and only if there exists $z^*\in\mathbb{R}^n$ satisfying $x^*=\mathrm{prox}_{\varphi}^{\Lambda}(z^*)$ and $z^*$ is a solution of the nonsmooth equation $F_{\mathrm{nor}}^{\Lambda}(z)=0$.
\end{opti con noraml map}

By Lemma~\ref{lem: partial psi}, the calculation of $\psi'(x;d_s(x))d_s(x)$ is equivalent to solving an optimization problem $\psi'(x;d_s(x))d_s(x)=\mathbf{P}_{\partial\psi(x)}(0)=\mathop{\argmin}_{v\in\partial \psi(x)} \Vert v\Vert$. Alternatively, we can first solve
\begin{equation}
    \tau(x)=\mathop{\argmin}_{z\in\mathbb{R}^n}\Vert F_{\mathrm{nor}}^{\Lambda}(z)\Vert \quad \text{s.t.}\quad \mathrm{prox}_{\varphi}^{\Lambda}(z)=x\label{def tau}
\end{equation}
and then compute $\psi'(x;d_s(x))d_s(x)=F_{\mathrm{nor}}^{\Lambda}(\tau(x))$. By the definition of $F_{\mathrm{nor}}^{\Lambda}(z)$, solving \eqref{def tau} is equivalent to
\begin{equation*}
    \tau(x)=\mathop{\argmin}_{z\in\mathbb{R}^n}\Vert \nabla f(x)+\Lambda(z-x)\Vert \quad \text{s.t.}\quad \mathrm{prox}_{\varphi}^{\Lambda}(z)=x,
\end{equation*}
which combined with \eqref{prox opti con} leads to
\begin{equation}
    \tau(x)=x+\Lambda^{-1}\mathbf{P}_{\partial \varphi(x)}(-\nabla f(x)).\label{comp tau}
\end{equation}
and
\begin{equation}
    \psi'(x;d_s(x))d_s(x)=F_{\mathrm{nor}}^{\Lambda}(\tau(x))=\nabla f(x)+\mathbf{P}_{\partial \varphi(x)}(-\nabla f(x)).\label{comp F tau}
\end{equation}
A closed form representation of the mapping $F_{\mathrm{nor}}^{\Lambda}(\tau(x))$ can be derived for $\ell_1$-optimization, group lasso, and $\ell_{\infty}$-optimization. We present $F_{\mathrm{nor}}^{\Lambda}(\tau(x))$ for an $\ell_1$-problem in Example~\ref{example-tau-L1}; other examples are summarized in the appendix in Example~\ref{example-tau}.

\newtheorem{example-tau-L1}[lem: partial psi]{Example}
\begin{example-tau-L1}[$F_{\mathrm{nor}}^{\Lambda}(\tau(x))$ for $\ell_1$-optimization]\label{example-tau-L1}
Suppose that $\varphi(x)=\Vert x\Vert_1$ and $\Lambda=\mathrm{diag}(\lambda_1,\lambda_2,\cdots,\lambda_n)$, then by \eqref{comp F tau}, we can compute
\begin{equation*}
F_{\mathrm{nor}}^{\Lambda}(\tau(x))_i=
\begin{cases}
\nabla f(x)_i-\mathbf{P}_{[-1,1]}(\nabla f(x)_i), &x_i=0,\\
\nabla f(x)_i+\mathrm{sgn}(x_i), &x_i\neq 0,\\
\end{cases}\quad \forall~i=1,2,\cdots,n.
\end{equation*}
\end{example-tau-L1}

\subsubsection{Natural Residual}
Another possible choice for $g(x)=u(x)d(x)$ can be based on the so-called \emph{natural residual},
 \begin{equation}
     F_{\mathrm{nat}}^{\Lambda}(x):=x-\mathrm{prox}_{\varphi}^{\Lambda}(x-\Lambda^{-1}\nabla f(x)).
 \end{equation}
 Similar to the normal map, $F_{\mathrm{nat}}^{\Lambda}$ can be used as a criticality measure.
 
 \newtheorem{opti con natural res}[lem: partial psi]{Lemma}
\begin{opti con natural res}
A point $x^*$ is a stationary point of problem \eqref{obj fun:composite} if and only if $x^*$ is a solution of the nonsmooth equation $F_{\mathrm{nat}}^{\Lambda}(x)=0.$
\end{opti con natural res}

Following \cite[Proposition 4.2]{FukMin81}, the directional derivative at $x$ along $-F_{\mathrm{nat}}^{\Lambda}(x)$ satisfies $\psi'(x;-F_{\mathrm{nat}}^{\Lambda}(x))\leq -\Vert F_{\mathrm{nat}}^{\Lambda}(x)\Vert_{\Lambda}^2$. Thus, the direction $$d(x)=\begin{cases}
  -\frac{F_{\mathrm{nat}}^{\Lambda}(x)}{\Vert F_{\mathrm{nat}}^{\Lambda}(x)\Vert_\Lambda} & \text{if } 0\notin\partial\psi(x),\\
  0 & \text{if } 0\in\partial\psi(x),\end{cases}$$ is a descent direction with the directional derivative
\begin{equation*}
\psi'(x;d(x))\leq -\frac{\Vert F_{\mathrm{nat}}^{\Lambda}(x)\Vert_{\Lambda}^2}{\Vert F_{\mathrm{nat}}^{\Lambda}(x)\Vert_\Lambda}=-\Vert F_{\mathrm{nat}}^{\Lambda}(x)\Vert_\Lambda.
\end{equation*}
We can choose $u(x)=-\Vert F_{\mathrm{nat}}^{\Lambda}(x)\Vert_\Lambda$, which implies that
\begin{equation}
    g(x)= F_{\mathrm{nat}}^{\Lambda}(x).\label{g(x):2}
\end{equation}

\subsection{Stepsize Safeguard}
If we utilize a smooth model $m_k$, once we have selected the descent direction $d$, (directionally) noncontinuous points of $\psi'(\cdot,d)$ will contribute to the inaccuracy of $m_k$. Hence, we should keep the stepsize relatively small to avoid those points. For any $x,d\in\mathbb{R}^n$ with $\Vert d\Vert=1$, we set the stepsize safeguard $\Gamma(x,d)$ to guarantee that $t \mapsto \psi'(x+td;d)$ is continuous on $(0,\Gamma(x,d))$, which   is equivalent to saying that $t \mapsto \psi\left(x+td\right)$ is continuously differentiable on $t\in (0, \Gamma(x,d))$.

We prefer to choose the largest possible value of $\Gamma(x,d)$:
\begin{equation}
\Gamma(x,d)=\Gamma_{\max}(x,d):=\sup\left\{T>0 :  \tilde{\psi}^\prime_{x,d}(t):=\psi^\prime(x+td;d)\in C^1(0, T)\right\},\label{large Gamma}
\end{equation}
since it can intuitively lead to faster convergence. We will see that this choice works well for polyhedral problems, where $\varphi$ is the supremum of several affine functions, such as in $\ell_1$- and $\ell_\infty$-optimization. 

However, in some other cases, we may need to set $\Gamma(x,d)$ more carefully. For example, for the group lasso problem $\min_{X \in \R^{n_1\times n_2}} f(X) + \varphi(X)$, 
where $f$ is smooth and $\varphi$ is given by $\varphi(X)=\sum_{i=1}^{n_2}\left\Vert X_i\right\Vert$, an appropriate choice for $X=(X_1,X_2,\cdots,X_{n_2})\in \mathbb{R}^{n_1\times n_2}$ and $D=(D_1,D_2,\cdots,D_{n_2})\in \mathbb{R}^{n_1\times n_2}$ with $\Vert D\Vert=1$ is:
\begin{equation}
    \Gamma(X,D)=\min\left\{\Gamma_{\max}(X,D),\min_{ X_i\neq 0}\frac{\Vert X_i\Vert^{1+\sigma}}{1-\theta_i^2},\min_{X_i\neq 0}\frac{\Vert X_i\Vert}{\max\{-2\theta_i,0\}}\right\}, \quad \sigma > 0.\label{gp lasso Gamma}
\end{equation}
Here, $\theta_i$ is given by $\theta_i:=\langle X_i,D_i\rangle /(\Vert X_i\Vert \cdot\Vert D_i\Vert)$ and we use $c/0:=+\infty$ if $c>0$. The term $\Gamma_{\max}(X,D)$ is defined as in \eqref{large Gamma}. This $\Gamma(X,D)$ is specifically designed to overcome some technical difficulties; see, e.g., Lemma~\ref{verify:B.3}. 

Next, let us define the function $\Gamma:\mathbb{R}^n\rightarrow\mathbb{R}^+$ via
\begin{equation}
    \Gamma(x):=\inf_{d\in\mathbb{R}^n,\ \Vert d\Vert=1}\Gamma_{\max}(x,d).\label{Gamma(x)}
\end{equation}
The scalar $\Gamma(x)$ is important in our convergence analysis as it provides a lower bound for the stepsize safeguard. For the composite program \eqref{obj fun:composite} with $\varphi(x)=\Vert x\Vert_1$, $\Gamma$ can be simply calculated as follows $\Gamma(x)=\min\{|x_i|: x_i\neq 0\}$ where $\min\emptyset:=+\infty$. Further examples for $\Gamma$ can be found in Appendix~\ref{example-gamma}.

\subsection{Truncation Operators}
\label{sec:truncatable function}
Since in our algorithmic design we utilize simple, linear-quadratic models to approximate the nonsmooth function $\psi$, we need to introduce stepsize safeguards that allow to intrinsically control the accuracy of the model. However, if the ``safeguard'' $\Gamma(x)$ is very small, the resulting step might be close to the old iterate and the algorithm can start to stagnate. In order to prevent such an undesirable behavior, we discuss an additional modification step that allows to increase $\Gamma(x)$. 

Specifically, given a point $x$, first we want to find a point $x'$ near $x$ such that $\Gamma(x')$ is relatively large. Let us consider the simplest case where $\psi=f+\varphi$ and $\varphi(x)=\Vert x\Vert_1$. If $x$ has a nonzero component with small absolute values, then $\Gamma(x)$ is also small. So we can replace those components with 0 and get a new point $x'$ satisfying $\Gamma(x')>\Gamma(x)$. Since only some components with small absolute values are truncated to 0, the point $x'$ is close to $x$. In more general cases, we define a class of functions that allow similar operations:

\newtheorem{def: truncation}[lem: partial psi]{Definition} 
\begin{def: truncation} \label{def: truncation}
Suppose that there exist a finite sequence $\{S_i\}_{i=0}^m$ satisfying $\mathbb{R}^n=S_0\supset S_1\cdots\supset S_m$, $\delta\in (0,+\infty]$, $\kappa>0$, and a function $T:\mathbb{R}^n\times(0,\delta]\rightarrow\mathbb{R}^n$ with following properties:

\begin{itemize}
    \item[{(i)}] $\Gamma(x)\geq \delta,\ \forall\  x\in S_m$;
    \item[{(ii)}] For any $a\in(0,\delta]$ and $x\in S_i\backslash S_{i+1}$, $i\in\{0,1,\cdots,m-1\}$, if $\Gamma(x)\geq a$, it holds that $T(x,a)=x$, otherwise we have $T(x,a)\in S_{i+1}$, $\Gamma(T(x,a))\geq a$, and $\Vert T(x,a)-x \Vert\leq \kappa a$.
\end{itemize}
Then we say that $\psi$ can be truncated and that $T$ is a truncation operator.
\end{def: truncation}

In Definition~\ref{def: truncation}, $\Gamma(T(x,a))\geq a$ means that we can make the value of $\Gamma(\cdot)$ larger by performing truncation and $\Vert T(x,a)-x \Vert\leq \kappa a$ implies that the change caused by $T(\cdot,a)$ can be controlled. 
Example~\ref{example-truncation-L1} shows that $\varphi(x)=\Vert x\Vert_1$ can be truncated and we present more examples ($\ell_\infty$-optimization and group lasso) in the appendix.

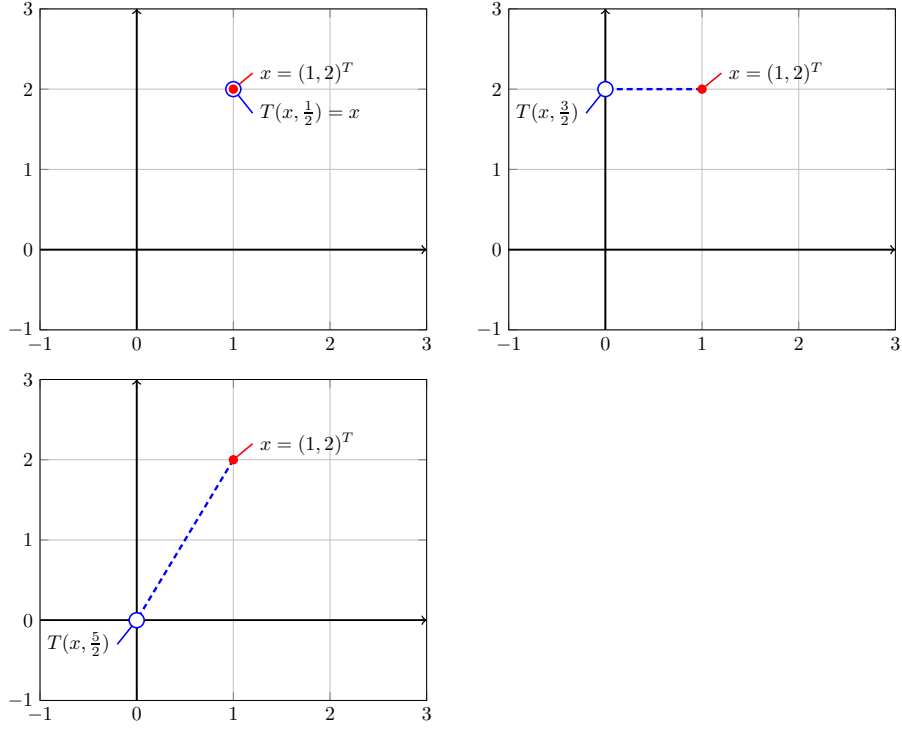
\begin{figure}[t]
\centering
\setlength{\belowcaptionskip}{-6pt}
\begin{tabular}{cc}
\begin{tikzpicture}[scale=0.75]
\begin{axis}[samples = 50, xmin = -1, xmax = 3, ymin = -1, ymax = 3,grid]
  \addplot[black,thick] coordinates { (-1,0) (3,0) };
  \addplot[black,thick] coordinates { (0,-1) (0,3) };
  \draw[black,thick,->] (axis cs:-1,0) -- (axis cs:3,0);
  \draw[black,thick,->] (axis cs:0,-1) -- (axis cs:0,3);
  \draw[red,thick] (axis cs:1,2) -- (axis cs:1.2,2.2) node[right,black]{$x = (1,2)^T$};
  \draw[blue,thick] (axis cs:1,2) -- (axis cs:1.2,1.7) node[right,black]{$T(x,\frac12) = x$};
  \node[circle,draw=blue,fill=white,scale=0.8,thick] at (axis cs:1,2) {};
  \node[circle,fill=red,scale=0.5] at (axis cs:1,2) {};
\end{axis}
\end{tikzpicture} &
\begin{tikzpicture}[scale=0.75]
\begin{axis}[samples = 50, xmin = -1, xmax = 3, ymin = -1, ymax = 3,grid]
  \addplot[black,thick] coordinates { (-1,0) (3,0) };
  \addplot[black,thick] coordinates { (0,-1) (0,3) };
  \draw[black,thick,->] (axis cs:-1,0) -- (axis cs:3,0);
  \draw[black,thick,->] (axis cs:0,-1) -- (axis cs:0,3);
  \draw[red,thick] (axis cs:1,2) -- (axis cs:1.2,2.2) node[right,black]{$x = (1,2)^T$};
  \draw[blue,thick] (axis cs:0,2) -- (axis cs:-0.2,1.7) node[left,black]{$T(x,\frac32)$};
  \draw[blue,very thick,densely dashed] (axis cs:1,2) -- (axis cs:0,2);
  \node[circle,draw=blue,fill=white,scale=0.8,thick] at (axis cs:0,2) {};
  \node[circle,fill=red,scale=0.5] at (axis cs:1,2) {};
\end{axis} 
\end{tikzpicture} \\
\begin{tikzpicture}[scale=0.75]
\begin{axis}[samples = 50, xmin = -1, xmax = 3, ymin = -1, ymax = 3,grid]
  \addplot[black,thick] coordinates { (-1,0) (3,0) };
  \addplot[black,thick] coordinates { (0,-1) (0,3) };
  \draw[black,thick,->] (axis cs:-1,0) -- (axis cs:3,0);
  \draw[black,thick,->] (axis cs:0,-1) -- (axis cs:0,3);
  \draw[red,thick] (axis cs:1,2) -- (axis cs:1.2,2.2) node[right,black]{$x = (1,2)^T$};
  \draw[blue,thick] (axis cs:0,0) -- (axis cs:-0.2,-0.3) node[left,black]{$T(x,\frac52)$};
  \draw[blue,very thick,densely dashed] (axis cs:1,2) -- (axis cs:0,0);
  \node[circle,draw=blue,fill=white,scale=0.8,thick] at (axis cs:0,0) {};
  \node[circle,fill=red,scale=0.5] at (axis cs:1,2) {};
\end{axis} 
\end{tikzpicture}
\end{tabular}
\caption{Illustration of the truncation operator in Example \ref{example-truncation-L1} for $\varphi(x)=\Vert x\Vert_1$ and $n=2$.}
\label{fig: truncation-L1}
\end{figure}

\newtheorem{example-truncation-L1}[lem: partial psi]{Example}
\begin{example-truncation-L1}[$\varphi(x)=\Vert x\Vert_1$]\label{example-truncation-L1}
For $i=0,1,\cdots,n$, we set $S_i=\{x\in\mathbb{R}^n\mid \mathrm{card}\{j=1,2,\cdots,n\mid x_j=0\}\geq i\}$, $m=n$, $\delta=+\infty$, $\kappa=\sqrt{n}$, and
\begin{equation*}
    T(x,a)_j=\mathbbm{1}_{|\cdot|\geq a}(x_j)x_j,\ j=1,2,\cdots,n,
\end{equation*}
where $\mathbbm{1}_A(\cdot)$ is the indicator function. Figure \ref{fig: truncation-L1} shows the truncation operator for $\varphi(x)=\Vert x\Vert_1$ and $n=2$ explicitly, where $S_1=\{(x_1,x_2)\mid x_1 x_2=0\}$ and $S_2=\{(0,0)\}$.
\end{example-truncation-L1}

Let us mention that, for a smooth regularizer $\varphi$, all properties discussed above are satisfied, since the stepsize safeguard can be chosen as $+\infty$ and no truncation is needed.

\section{A Nonsmooth Trust-region Method}
  \label{sec: alg}
  In this section, we present the algorithmic framework of our trust-region type method. The traditional trust-region framework using the pseudo-gradient $g(x)$ is employed in our algorithm with potential refinement using a stepsize safeguard. The traditional framework is standard but it might not be accurate enough to pass the acceptance test. In each iteration, we first perform a classical trust-region step. A stepsize safeguard for refinement is used if the traditional step fails. In order to promote large stepsizes, a novel truncation step is proposed.
  
  We first introduce the model function and trust-region subproblem. Then, we propose several modification steps including the choice of the stepsize and the novel truncation step. The final algorithm is presented at the end of this section. We also present some methods for solving the corresponding trust-region subproblem in Appendix~\ref{solution_sub}.

\subsection{Model Function and Trust-region Subproblem}
Recall that in the classical trust-region method for a smooth optimization problem, $\min_{x\in \mathbb{R}^n} \psi(x)$, the model function is $m_k(s)=\psi(x^k)+\langle\nabla \psi(x^k),s\rangle+\frac{1}{2}\langle s,B^k s\rangle$. As mentioned at the beginning of Section~\ref{sec: descent directions}, a natural extension is
\begin{equation}
    m_k(s)=\psi(x^k)+\langle \psi'(x^k;d_s(x^k))d_s(x^k),s\rangle+\frac{1}{2}\langle s ,B^k s\rangle.
\end{equation}
This model function is still quadratic and fits the objective function well along the steepest descent direction, which means that they have the same directional derivative in this direction. Though approximating the nonsmooth function $\psi$ with a quadratic function might not lead to good trust-region models in general, we can design specific quadratic models
that fit $\psi$ well along certain directions.

One may wish to use different descent directions. Or it might be expensive to compute the steepest descent direction $d_s(x)$ and its Clarke's generalized directional derivative $\psi'(x;d_s(x))$. Therefore we use a descent direction $d(x)$ satisfying \eqref{d(x)} and $u(x)$ satisfying \eqref{u(x)} instead. We can now define our model function
\begin{equation}
    m_k(s)=\psi_k+\langle g^k,s\rangle+\frac{1}{2}\langle s,B^k s\rangle, \label{eq_mk}
\end{equation}
where $\psi_k=\psi(x^k)$ and $g^k=g(x^k)=u(x^k)d(x^k)$, and the associated trust-region subproblem is given by
\begin{equation}
\min_s\ m_k(s)=\psi_k+\langle g^k,s\rangle+\frac{1}{2}\langle s,B^k s\rangle \quad \text{s.t.} \quad \Vert s\Vert\leq\Delta_k \label{sub-problem}.
\end{equation}
This subproblem is quadratic and coincides with the classical approaches if $B^k$ is symmetric. The matrices $B^k$ in the model \eqref{eq_mk} are typically chosen to capture or approximate the second-order information of the objective function $\psi$. However, such a careful choice is only required when certain local convergence properties should be guaranteed. Similar to classical trust-region methods, global convergence of the method will generally not be affected by $\{B^k\}$ and flexible choices of $B^k$ are possible under some mild boundedness conditions that will be specified later.

We remark that the major difference between our algorithm and other nonsmooth trust-region methods in the literature is the first-order information in the model function. The methods in \cite{Nonsmooth-TR-Qi-Sun-94, Nonsmooth-TR-Christof-18} employ first-order terms that tend to be more complicated than a simple linear function, in order to approximate the nonsmooth objective function well and to satisfy certain accuracy assumptions. Although the model function in \cite{TR-quadratic-model-Akbari-15} is quadratic, their first-order term needs to be built using a steepest descent direction based on the $\epsilon$-subdifferential of $\psi$.

An important concept for solving \eqref{sub-problem} is the so-called \emph{Cauchy point}, which is defined via
\begin{equation*}
    s^k_C:=-\alpha_k^C g^k\ \text{and}\ \alpha_k^C:= \mathop{\argmin}_{0\leq t\leq \Delta_k/\Vert g^k\Vert} m_k(-t g^k).
\end{equation*}
The Cauchy point is computational inexpensive \cite[Algorithm 4.2]{Num-Opti-Nocedal-06} and it leads to sufficient reduction of the model function (Cauchy decrease condition):
\begin{equation}
    m_k(0)-m_k(s^k_C)\geq \frac{1}{2}\Vert g^k\Vert\min\left\{\Delta_k, \frac{\Vert g^k\Vert}{\Vert B^k\Vert}\right\},\label{Cau_Reduce}
\end{equation}
see, e.g., in \cite[Lemma 4.3]{Num-Opti-Nocedal-06}. Furthermore, it can be shown that
\begin{equation}
    m_k(0)-m_k(\alpha_k\bar{s}^k_C)\geq \frac{\alpha_k}{\Vert s^k_C\Vert}\cdot\frac{1}{2}\Vert g^k\Vert\min\left\{\Delta_k, \frac{\Vert g^k\Vert}{\Vert B^k\Vert}\right\},\label{Cau_Reduce_step}
\end{equation}
where $\bar{s}^k_C=s^k_C/\Vert s^k_C\Vert$ and $0<\alpha_k\leq \Vert s^k_C\Vert$.

In our algorithm, we need to generate an approximate solution of \eqref{sub-problem} that achieves a similar model descent compared to the Cauchy descent condition \eqref{Cau_Reduce} in some sense. More precisely, we need to recover a solution $s^k$ satisfying
\begin{equation}
    m_k(0)-m_k(s^k)\geq \frac{\gamma_1}{2}\Vert g^k\Vert\min\left\{\Delta_k, \gamma_2\Vert g^k\Vert\right\},\label{model reduce0}
\end{equation}
where $\gamma_1,\gamma_2>0$ are constants which do not depend on $k$, and
\begin{equation}
    m_k(0)-m_k(s^k)\geq (1-\ell(\Vert s^k\Vert))[ m_k(0)-m_k(s^k_C) ],\label{model reduce}
\end{equation}
where $\ell:\mathbb{R}^+\rightarrow[0,1]$ is chosen as a monotonically decreasing function with $\lim_{\Delta\rightarrow 0^+}\ell(\Delta)=0$. The classical choice $\ell(\Delta)\equiv 0$ is also allowed.

\subsection{Suitable Stepsizes}
In the trust-region framework, we will work with the parameters $0<\eta \leq \eta_1<\eta_2<1$, $0<r_1<1<r_2$, and $\Delta_{\max}>0$. Let $s^k$ denote the generated solution of \eqref{sub-problem}.
Similar to the classical trust-region method, we define the ratio between actual reduction and predicted reduction as
\begin{equation}
\rho^1_k=\frac{\psi(x^k)-\psi(x^k+s^k)}{m_k(0)-m_k(s^k)}.\label{rho-1}
\end{equation}
If the proposed step $x^k+s^k$ is ``successful'', i.e., $\rho^1_k\geq \eta_1$, we accept the step, i.e., $\tilde{x}^k=x^k+s^k$, and update the trust-region radius $\Delta_k$ as
\begin{equation}
\Delta_{k+1}=\begin{cases}\min\{{\Delta}_{\max},r_2\Delta_k\} & \text{if } \rho^1_k>\eta_2,\\
\Delta_k & \text{otherwise}.\end{cases}\label{Update:Delta-1}
\end{equation}
If $\rho^1_k< \eta_1$, we can introduce an additional stepsize strategy to refine the step. Specifically, we now consider the normalized descent direction $\bar{s}^k:={s^k}/{\Vert s^k\Vert}$. In the following, we will always use the notation $\bar{x} := {x}/{\Vert x\Vert}$. Instead of setting the stepsize as $\Vert s^k\Vert$ and working with $s^k$ directly, we calculate $\alpha_k$ via
\begin{equation}
    \alpha_k=\min\left\{\Gamma(x^k,\bar{s}^k),\Vert s^k\Vert\right\},\label{stepsize}
\end{equation}
where $\Gamma(x^k,\bar{s}^k)$ is the stepsize safeguard. If
\begin{equation}
m_k(0)-m_k(\alpha_k\bar{s}^k)\geq \frac{\alpha_k}{2\Vert s^k\Vert}(m_k(0)-m_k(s^k)),\label{model_convex}
\end{equation}
which means that the modified step yields sufficient descent, we use the direction $\bar{s}^k$ and the stepsize $\alpha_k$; otherwise we set 
\begin{equation}
s^k=s^k_C,\ \bar{s}^k=\bar{s}^k_C,\text{ and }\alpha_k=\min\left\{\Gamma(x^k,\bar{s}^k_C),\Vert s^k_C\Vert\right\}.\label{stepsize2}
\end{equation}
If $m_k$ is convex (which, e.g., can be ensured when the matrix $B^k$ is chosen to be positive semidefinite) then \eqref{model_convex} holds automatically and the latter case will not occur. In \eqref{stepsize2} we utilize the Cauchy point and the corresponding stepsize as a simpler gradient-based step. As we have seen such a step can always guarantee certain descent properties. Next, we perform a second ratio test
\begin{equation}
\rho^2_k=\frac{\psi(x^k)-\psi(x^k+\alpha_k\bar{s}^k)}{m_k(0)-m_k(\alpha_k\bar{s}^k)}.\label{rho-2}
\end{equation}
According to this ratio, we update the trust-region radius as
\begin{equation}
\Delta_{k+1}=\begin{cases}r_1\Delta_k & \text{if } \rho^2_k<\eta_1,\\
\min\{{\Delta}_{\max},r_2\Delta_k\} & \text{if } \rho^2_k>\eta_2,\\
\Delta_k & \text{otherwise},\end{cases}\label{Update:Delta-2}
\end{equation}
and decide whether to accept the proposed step
\begin{equation}
\tilde{x}^k=\begin{cases} x^k+\alpha_k\bar{s}^k & \text{if } \rho_k^2 \geq \eta,\\
x^k & \text{if } \rho_k^2 < \eta.\end{cases}\label{Update:x}
\end{equation}
We declare the step as ``subsuccessful'' if $\rho^1_k<\eta_1$ while $\rho^2_k\geq\eta$, i.e., even if the original step is unsuccessful, the refined version can still provide some descent which is essential to guarantee convergence.

\subsection{Truncation Step}

It might not be suitable to simply set $x^{k+1}=\tilde{x}^k$, since $\Gamma(\tilde{x}^k)$ can be very small and larger $\Gamma$-values increase the stepsize and improve the fitness of the model. Our idea is to allow a small modification of $\tilde{x}^k$ and to get a new point $x^{k+1}$ with relatively large $\Gamma(x_{k+1})$ although such modification may cause an increase of the objective function. In the following, we describe an algorithmic procedure for increasing the safeguard $\Gamma(x^{k+1})$ for functions which can be truncated.

Suppose that $\varphi$ can be truncated and let $S_0,S_1\cdots,S_m$, $\delta$, and $T$ be the corresponding truncation parameters and operators, respectively.
Let $\{\epsilon_s\}_{s=0}^{\infty} \in \ell_1^+$ be a positive and strictly decreasing sequence that is upper-bounded by $\delta$ as well as summable.  Since the sets $\{S_j: j=0,1,\cdots,m\}$ are nested and cover the whole $\mathbb{R}^n$, we know that there exists a unique index $i\in\{0,1,\cdots,m\}$  with $\tilde{x}^k \in S_i\backslash S_{i+1}$, where $S_{m+1}=\emptyset$. In the following, we define $N_j := S_j \backslash S_{j+1}$ and introduce a global counter $c_j$ that is associated with each set $N_j$ and that counts the total number of truncations performed on points in the set $N_j$ for $j=0,1,\cdots,m$. Depending on the safeguard $\Gamma(\tilde{x}^k)$ we then decide whether $\tilde{x}^k$ should be truncated via applying the truncation operator or not. The whole process is given as follows: find $i \in\{0,1,\cdots,m\}$ such that $\tilde{x}^k \in N_i$; if $\Gamma(\tilde{x}^k) < \epsilon_{c_i}$, set $\tilde{x}^k\leftarrow T(\tilde{x}^k, \epsilon_{c_i})$, otherwise we keep $\tilde{x}^k$ unchanged; and update $c_i = c_i + 1$ if $\tilde{x}^k$ is updated. 

 This procedure is repeated until $\Gamma(\tilde{x}^k) \geq \epsilon_{c_i}$. Lemma~\ref{truncation-finite-stop}  implies that this algorithm is well-defined and terminates within a finite number of steps. We call the whole procedure a \emph{truncation step} which is presented in Algorithm \ref{alg:truncation}.
 
 \begin{algorithm}[t]
 \caption{Truncation step}
\hspace*{0.02in} {\bf Input:} $\tilde{x}^k$ and $c_j,\ j=0,1,\cdots,m.$
\begin{algorithmic}[1]

\WHILE {true}
    \STATE Compute the unique $i$ such that $\tilde{x}^{k}\in S_i\backslash S_{i+1}$.
    \IF {$\Gamma(\tilde{x}^{k})<\epsilon_{c_i}$}
    \STATE $\tilde{x}^{k}\leftarrow T(\tilde{x}^{k}, \epsilon_{c_i})$.
    \STATE $c_i\leftarrow c_i+1$.
    \ELSE
        \STATE break.
    \ENDIF
\ENDWHILE

\end{algorithmic}
\hspace*{0.02in} {\bf Output:} $x^{k+1}=\tilde{x}^k$ and $c_j,\ j=0,1,\cdots,m.$
\label{alg:truncation}
\end{algorithm}

\newtheorem{truncation-finite-stop}[lem: partial psi]{Lemma}
\begin{truncation-finite-stop}
Algorithm \ref{alg:truncation} will terminate in at most $m$ steps.
\label{truncation-finite-stop}
\end{truncation-finite-stop}

\begin{proof}
Since for any $x\in S_m$ and $s\in \mathbb{N}$, we have $\Gamma(x)\geq \delta \geq \epsilon_s$ and the operator $T$ moves points in $S_i\backslash S_{i+1}$ into $S_{i+1}$, $T$ is performed on $\tilde{x}_k$ at most $m$ times before Algorithm \ref{alg:truncation} terminates.
\end{proof}

The iterate $\tilde{x}^k$ will be changed when performing Algorithm \ref{alg:truncation}. For simplicity, in the rest of this paper, when we mention $\tilde{x}^k$, we always mean the input of Algorithm \ref{alg:truncation}.

\subsection{Algorithmic Framework}
\label{sec:alg framework}
We now present a nonsmooth trust-region framework with quadratic model functions that combines the mentioned strategies. One of  the main advantages is that the corresponding subproblem can be cheaply formulated and solved. Specifically, the first-order term of our model can be constructed using any kind of descent direction and the second-order term $B^k$ is only required to satisfy a classical boundedness condition to guarantee global convergence. 
Moreover, the resulting trust-region subproblem coincides with the classical one and can be solved using classical methods. Let us further note that in order to obtain fast local convergence, $g^k$ and $B^k$ need to be chosen and coupled in a more careful way. This is explored in more detail in Section~\ref{sec: loc_conv} and \ref{numerics}.

The full algorithm is shown in Algorithm \ref{alg:main}. We require the following parameters: $0<\eta<\eta_1<\eta_2<1$, $0<r_1<1<r_2$, ${\Delta}_{\max} >0$, $\gamma_1>0$, $\gamma_2>0$, and a positive and strictly decreasing sequence $\{\epsilon_s\}_{s=0}^{\infty} \in \ell_1^+$ which is upper-bounded by $\delta$ in Definition~\ref{def: truncation}. We also assume that there is a monotonically decreasing function $\ell:\mathbb{R}^+\rightarrow[0,\frac{1}{2}]$ with $\lim_{\Delta\rightarrow 0^+}\ell(\Delta)=0$. Some additional discussions of how to solve the trust-region subproblems \eqref{sub-problem} can be found in the appendix.

\begin{algorithm}[t]
\caption{A trust-region method for nonsmooth nonconvex optimization}
{\bf Initialization:} initial point $x^0\in\mathbb{R}^n$, initial trust-region radius $\Delta^0$, iteration $k:=0$, global counters $c_j=0,\ j=0,1,\cdots,m$. 

\begin{algorithmic}[1]
\WHILE {$\|g^k\| \neq 0$}
    \STATE Compute $d(x^k)$, $u(x^k)$ and $g^k=u(x^k)d(x^k)$ and choose $B_k \in \R^{n \times n}$.
    \STATE Solve the trust-region subproblem \eqref{sub-problem} and obtain $s^k$ that satisfies \eqref{model reduce0} and \eqref{model reduce}.
    \STATE Compute $\rho^1_k$ according to \eqref{rho-1}.
    \IF{$\rho^1_k\geq \eta_1$}
        \STATE $\tilde{x}^k:=x^k+s^k$.
        \STATE Compute $\Delta_{k+1}$ according to \eqref{Update:Delta-1}.
    \ELSE
        \STATE Compute $s^k$, $\bar{s}^k$, and $\alpha_k$ according to \eqref{stepsize}, \eqref{model_convex}, and \eqref{stepsize2}.
        \STATE Compute $\rho^2_k$ according to \eqref{rho-2}.
        \STATE Compute $\Delta_{k+1}$ according to \eqref{Update:Delta-2}.
        \STATE Compute $\tilde{x}^{k}$ according to \eqref{Update:x}.
    \ENDIF
    \STATE Perform Algorithm \ref{alg:truncation}, get $x^{k+1}$ and update $c_j,\ j=0,1,\cdots,m$.
    \STATE $k\leftarrow k+1$.
\ENDWHILE
\end{algorithmic}\label{alg:main}
\end{algorithm}

We want to mention here that, for iteration $k\geq 1$, if $x^k+\alpha_k\bar{s}^k$ is not accepted, i.e., $\tilde{x}^k=x^k$, we have $x^{k+1}=\tilde{x}^k=x^k$, which means that no truncation is performed on $\tilde{x}^k$. This is because $x^k$ satisfies the stopping criteria of Algorithm~\ref{alg:truncation} since it was the output of Algorithm~\ref{alg:truncation} in the last iteration.

\section{Global Convergence}
\label{sec: glo_conv}
In this section, we show the global convergence of Algorithm~\ref{alg:main}. Specifically, we will prove that every accumulation point of the sequence generated by Algorithm~\ref{alg:main} is a stationary point under some suitable assumptions and that the natural residual converges to zero along the generated iterates. 

\subsection{Assumptions}
In this subsection, we state the assumptions required for the convergence. Assumption~\ref{Assumption:A} summarizes the conditions on the objective function $\psi$ and its pseudo-gradient $g$.

\newtheorem{Assumption:A}[lem: partial psi]{Assumption}
\begin{Assumption:A}\label{Assumption:A}
We assume that $\psi$ and $g$ have the following properties:

\begin{itemize}[leftmargin=7ex]
\item[{(A.1)}\label{A.1}] $\psi$ is bounded from below by $L_b$.
\item[{(A.2)}\label{A.2}] If $g(x)\neq 0$, then there exists $r,\epsilon>0$ such that $\Vert g(y)\Vert\geq \epsilon,\ \forall\  y\in B_r(x)$.
\end{itemize}
\end{Assumption:A}

Assumption (A.1) is a standard assumption. Assumption (A.2) means that the first-order model will not vanish sharply. Condition (A.2) holds automatically if $x \mapsto \Vert g(x)\Vert$ is lower semicontinuous.

\newtheorem{verify:A.2}[lem: partial psi]{Lemma}
\begin{verify:A.2} Assumption (A.2) is satisfied for the choices in \eqref{g(x):1} and \eqref{g(x):2}.\label{verify:A.2}
\end{verify:A.2}

\begin{proof}
(i) For \eqref{g(x):1}, set $\epsilon=\frac{1}{2}\Vert F_{\mathrm{nor}}^{\Lambda}(\tau(x))\Vert>0$. Suppose that there exist a sequence $\{y^m\}_{m}$ satisfying $y^m\rightarrow x$ and $\Vert F_{\mathrm{nor}}^{\Lambda} (\tau(y^m))\Vert<\epsilon$ for all $m\in\mathbb{N}$. By the local boundedness of $\partial\varphi$, we can infer that $\{\mathbf{P}_{\partial\varphi(y^m)}(-\nabla f(y^m))\}_{m}$ is bounded and thus, $\{\mathbf{P}_{\partial\varphi(y^m)}(-\nabla f(y^m))\}_{m=0}^{\infty}$ has a convergent subsequence. Without loss of generality, we assume that the whole sequence $\{\mathbf{P}_{\partial\varphi(y^m)}(-\nabla f(y^m))\}_{m=0}^{\infty}$ converges. Let us set $w=\lim_{m\rightarrow\infty}\mathbf{P}_{\partial\varphi(y^m)}(-\nabla f(y^m))$. Using the upper semicontinuity or the closedness of $\partial\varphi$, $y^m\rightarrow x$, and $\partial\varphi(y^m)\ni \mathbf{P}_{\partial\varphi(y^m)}(-\nabla f(y^m))\rightarrow w$, it follows $w\in\partial\varphi(x)$. Therefore, we obtain $\nabla f(x)+w\in\partial\psi(x)$ with $\Vert \nabla f(x)+w\Vert\leq \epsilon<\Vert F_{\mathrm{nor}}^{\Lambda}(\tau(x))\Vert$, which contradicts the optimality of $F_{\mathrm{nor}}^{\Lambda}(\tau(x))$. We can conclude that for some $r>0$, it holds $\Vert g(y)\Vert=\Vert F_{\mathrm{nor}}^{\Lambda} (\tau(y))\Vert=-\psi'(y,d_s(y))\geq \epsilon$ for all $y\in B_{r}(x)$. Hence, condition (A.2) is satisfied.

(ii) For \eqref{g(x):2}, assumption (A.2) is a simple consequence of the continuity of $F^\Lambda_{\mathrm{nat}}$. 
\end{proof}

Next, we present several assumptions on the iterates and sequences generated by Algorithm~\ref{alg:main}. 

\newtheorem{Assumption:B}[lem: partial psi]{Assumption}
\begin{Assumption:B}
Let $\{x^k\}$ and $\{B^k\}$ be generated by Algorithm~\ref{alg:main}. We assume:
\begin{itemize}[leftmargin=7ex]
    \item[{(B.1)}\label{B.1}] $\{x^k\}$ is bounded, i.e., there exist $R>0$ with $\{x^k\}\subseteq B_R(0)$.
    \item[{(B.2)}\label{B.2}] There exists $\kappa_B>0$ with $\sup_{k\in\mathbb{N}}\Vert B^k\Vert\leq\kappa_B<\infty$.
    \item[{(B.3)}\label{B.3}] For any subsequence $\{k_\ell\}_{\ell=0}^{\infty}\subseteq\mathbb{N}$, if $\{x^{k_\ell}\}$ converges and we have $\alpha_{k_\ell}\rightarrow 0$, then it holds that
\end{itemize}
    \begin{equation}
        \psi(x^{k_\ell}+\alpha_{k_\ell} \bar{s}^{k_\ell})-\psi(x^{k_\ell})-\alpha_{k_\ell} \psi'(x^{k_\ell};\bar{s}^{k_\ell}) = o( \alpha_{k_\ell}) \quad \ell \to \infty.\label{B.3 psi}
    \end{equation}
\begin{itemize}[leftmargin=7ex]
    \item[{(B.4)}\label{B.4}] For every $\epsilon>0$ there is $\epsilon'>0$ such that for all $x^k$ with $\Gamma(x^k) \geq \epsilon$ it follows $\Gamma\left(x^k,\bar{s}^k\right)\geq \epsilon'$.
\end{itemize}
\end{Assumption:B}

The conditions (B.1)--(B.3) are standard assumptions. Condition (B.2) is frequently used in classical trust-region theory, see, e.g., \cite[Theorem 4.5]{Num-Opti-Nocedal-06}. Assumption (B.3) is required to ensure uniform accuracy of the model function. Similar assumptions also appear in other convergence analyses of nonsmooth trust-region methods. For instance, condition A.2 in \citep{Nonsmooth-TR-Qi-Sun-94} and assumption 2.4(2c) in \cite{Nonsmooth-TR-Christof-18} have similar formats. However, as far as we can tell, a simple modification of the convergence analysis in previous works (such as \cite{Nonsmooth-TR-Qi-Sun-94, Nonsmooth-anal-opti-Clason-18}) may not prove the convergence of our method. This is partly because assumption (B.3) is only required when the stepsize $\alpha_k$ is smaller than or equal to some safeguard, see also \eqref{stepsize}. Moreover, even when the trust-region radius is relatively large and we pass the acceptance, the potential descent could be small. This is also further motivates our truncation step which allows to enlarge the stepsize.

Let us notice that condition (B.3) is similar to but weaker than the concept of ``strict models'' introduced by Noll in \citep{Cutting-Plane-Oracles-Noll-10}. In particular, assumption (B.3) only needs to hold at accumulation points of $\{x^k\}$ and along the specific and associated directions $\{\bar s^k\}$ while typical strict model assumptions are formulated uniformly for all points and directions in $\mathbb{R}^n$. In the following lemmas, we verify the  conditions (B.3) and (B.4) for two exemplary cases.

\newtheorem{verify:B.3}[lem: partial psi]{Lemma} 
\begin{verify:B.3} Suppose that $\psi$ is given by $\psi=f+\varphi$ with $\nabla f$ being locally Lipschitz continuous and assume that condition (B.1) holds. Then, (B.3) is satisfied in the following cases:
\begin{itemize}
\item[(i)] The mapping $\varphi$ is polyhedral and we have $\Gamma(x,d)\leq \Gamma_{\max}(x,d)$.
\item[(ii)] The problem is in the group lasso format and we set $\Gamma(x,d)$ as in \eqref{gp lasso Gamma}.
\end{itemize}
\label{verify:B.3}
\end{verify:B.3}

\begin{proof}
The local Lipschitz continuity of $\nabla f$ and the boundedness of $\{x^k\}$ imply that for any subsequence $\{k_\ell\}_{\ell=0}^{\infty}\subseteq\mathbb{N}$ with $x^{k_\ell}\rightarrow x$ and $\alpha_{k_\ell}\rightarrow 0$, it holds that $f(x^{k_\ell}+\alpha_{k_\ell} \bar{s}^{k_\ell})-f(x^{k_\ell})-\alpha_{k_\ell} f'(x^{k_\ell};\bar{s}^{k_\ell}) = o( \alpha_{k_\ell})$. Thus, it suffices to prove
\begin{equation}
\varphi(x^{k_\ell}+\alpha_{k_\ell} \bar{s}^{k_\ell})-\varphi(x^{k_\ell})-\alpha_{k_\ell} \varphi'(x^{k_\ell};\bar{s}^{k_\ell}) = o( \alpha_{k_\ell}).\label{B.3:varphi}
\end{equation}
If $\varphi$ is polyhedral, the function $\tilde{\varphi}_{x,d}(t):=\varphi\left(x+td\right)$ is linear on $(0,\Gamma_{\max}(x,d))$. Thus, \eqref{B.3:varphi} holds with the right side of the equality taken as zero.
    
In case of the group lasso problem and using the definition \eqref{gp lasso Gamma}, we can see that
$\Vert X_i^{k_\ell}\Vert+\alpha_{k_\ell} \theta_i^{k_\ell}\geq 0.5 \Vert X_i^{k_\ell}\Vert$, where $\theta_i^{k_\ell}:= \langle X_i^{k_\ell}, S_i^{k_\ell}\rangle  / (\Vert X_i^{k_\ell}\Vert \cdot\Vert S_i^{k_\ell}\Vert)$. Thus, we have
\begin{align*}
    \varphi(X^{k_\ell}+\alpha_{k_\ell} \bar{S}^{k_\ell})&-\varphi(X^{k_\ell})-\alpha_{k_\ell} \varphi'(X^{k_\ell};\bar{S}^{k_\ell})\\
    =&\sum_{X_i^{k_\ell}\neq 0}\Vert X_i^{k_\ell}+\alpha_{k_\ell} \bar{S}^{k_\ell}_i\Vert-\Vert X_i^{k_\ell}\Vert-\alpha_{k_\ell} \langle \bar{X}_i^{k_\ell}, \bar{S}^{k_\ell}_i \rangle \\
    =&\sum_{X_i^{k_\ell}\neq 0}(\Vert X_i^{k_\ell}\Vert^2+\alpha_{k_\ell}^2+2\alpha_{k_\ell}\Vert X_i^{k_\ell}\Vert \theta_i^{k_\ell})^{1/2}-\Vert X_i^{k_\ell}\Vert-\alpha_{k_\ell} \theta_i^{k_\ell}\\
    =&\sum_{X_i^{k_\ell}\neq 0}((\Vert X_i^{k_\ell}\Vert+\alpha_{k_\ell} \theta_i^{k_\ell})^2+\alpha_{k_\ell}^2(1-(\theta_i^{k_\ell})^2))^{1/2}-\Vert X_i^{k_\ell}\Vert-\alpha_{k_\ell} \theta_i^{k_\ell}\\
    \leq&\sum_{X_i^{k_\ell}\neq 0}\frac{\alpha_{k_\ell}^2(1-(\theta_i^{k_\ell})^2)}{2(\Vert X_i^{k_\ell}\Vert+\alpha_{k_\ell} \theta_i^{k_\ell})}\leq \sum_{X_i^{k_\ell}\neq 0}\alpha_{k_\ell}^2\cdot\frac{1-(\theta_i^{k_\ell})^2}{\Vert X_i^{k_\ell}\Vert},
\end{align*}
where the penultimate inequality follows from $(a^2+b)^{1/2}-a\leq {b}/{(2a)}$ for $a>0$, $b>0$. For all $i$, if $\lim_{\ell\rightarrow\infty} X_i^{k_\ell}\neq 0$, we have $\alpha_{k_\ell}^2/ {\Vert X_i^{k_\ell}\Vert}=o(\alpha_{k_\ell})$. If $0\neq X_i^{k_\ell}\rightarrow 0$, definition \eqref{gp lasso Gamma} implies $\alpha_{k_\ell}^2(1-(\theta_i^{k_\ell})^2)/{\Vert X_i^{k_\ell}\Vert} \leq \alpha_{k_\ell}\Vert X_i^{k_\ell}\Vert^\sigma = o(\alpha_{k_\ell})$. Therefore, condition \eqref{B.3:varphi} is satisfied.
\end{proof}

\newtheorem{verify:B.4}[lem: partial psi]{Lemma} 
\begin{verify:B.4}\label{verify:B.4} Condition (B.4) is satisfied for the choices \eqref{large Gamma} and \eqref{gp lasso Gamma} (for group lasso problems).
\end{verify:B.4}

\begin{proof} Using \eqref{large Gamma}, we immediately obtain $\Gamma(x,d)\geq \Gamma(x)$ for all $x$ and $d$ with $\Vert d\Vert=1$. Hence, in this case, we can set $\epsilon^\prime = \epsilon$. 

For the group lasso problem and \eqref{gp lasso Gamma}, Example \ref{example-gamma} establishes $\Gamma(X)=\min\left\{\Vert X_i\Vert: X_i\neq 0\right\}$. Consequently, from $\Gamma(X)\geq \epsilon$, we can deduce
\begin{equation*}
    \Gamma(X,D)=\min\left\{\Gamma_{\max}(X,D),\min_{ X_i\neq 0}\frac{\Vert X_i\Vert^{1+\sigma}}{1-\theta_i^2},\min_{X_i\neq 0}\frac{\Vert X_i\Vert}{\max\{-2\theta_i,0\}}\right\}\geq\min\{\epsilon^{1+\sigma},\epsilon/2\},
\end{equation*}
where $\theta_i:= {\langle X_i,D_i\rangle}/{(\Vert X_i\Vert \cdot\Vert D_i\Vert)}$ for any $D$ with $\Vert D\Vert=1$.
\end{proof}

\subsection{Convergence Analysis}
In this subsection, we will prove that every accumulation point of $\{x^k\}_{k}$ is a stationary point of \eqref{obj fun:composite}. First, in Lemma~\ref{B.3:global}, we derive a global version of assumption (B.3) over the ball $\overline{B_R(0)}$.

\newtheorem{B.3:global}[lem: partial psi]{Lemma}
\begin{B.3:global} Suppose that (B.1) and (B.3) are satisfied and that Algorithm~\ref{alg:main} does not terminate within finitely many steps. Let $\{x^k\}_{k}$ be a sequence generated by Algorithm~\ref{alg:main}. Then, there exists a function $h:(0,\infty)\rightarrow[0,\infty]$ with $\lim_{\Delta\rightarrow0^+}h(\Delta)=0$ and
\begin{equation}
    \psi(x^k+\alpha_k \bar{s}^k)-\psi(x^k)-\alpha_k \psi'(x^k;\bar{s}^k)\leq h(\Delta)\alpha_k\label{eq-B.3:global},
\end{equation}
for all $k$ with $\Delta_k\leq \Delta$. \label{B.3:global}
\end{B.3:global}

\begin{proof} We set
\begin{equation*}
    h(\Delta)=\max\left\{\sup_{j:\ \Delta_j\leq \Delta}\frac{\psi(x^{j}+\alpha_{j} \bar{s}^{j})-\psi(x^{j})- \alpha_{j} \psi'(x^{j};\bar{s}^{j})}{\alpha_j},0\right\}.
\end{equation*}
From the definition, it directly follows $h(\Delta^1)\leq h(\Delta^2)$ for all $0<\Delta^1<\Delta^2<\infty$. Thus, it suffices to show that for every $\epsilon>0$ there exists $\Delta>0$ such that $h(\Delta)\leq\epsilon$. Let us assume that for some $\epsilon>0$ we have $h(\Delta)>\epsilon$ for all $\Delta>0$. Then, there exists a subsequence $\{k_\ell\}_{\ell=0}^{\infty}\subseteq\mathbb{N}$, such that $\Delta_{k_\ell}\rightarrow 0$ and
\begin{equation}
    \frac{\psi(x^{k_\ell}+\alpha_{k_\ell} \bar{s}^{k_\ell})-\psi(x^{k_\ell})- \alpha_{k_\ell} \psi'(x^{k_\ell};\bar{s}^{k_\ell})}{\alpha_{k_\ell}}\geq \epsilon \quad \forall~\ell \in \N.\label{B.3:contradict}
\end{equation}
Since $\{x^{k_\ell}\}_{\ell}$ is bounded, it has a convergent subsequence $\{x^{k_{\ell_m}}\}_{m}$. Due to $\Delta_{k_\ell}\rightarrow 0$ it follows $\alpha_{k_{\ell_m}}\rightarrow 0$. Therefore, \eqref{B.3 psi} and \eqref{B.3:contradict} yield a contradiction.
\end{proof}

Recall that the iterates $x^{k+1}$ result from a possible truncation of the trust-region steps $\tilde{x}^k$. We now prove that these truncation steps and the potential increase of the objective function $\psi$ can be controlled.

\newtheorem{truncation:finite move}[lem: partial psi]{Lemma}
\begin{truncation:finite move} 
Suppose that $\psi$ can be truncated and that assumption (B.1) holds. Then, we have $\sum_{k=0}^{\infty}\Vert x^{k+1}-\tilde{x}^k\Vert\leq m\kappa\sum_{i=0}^{\infty}\epsilon_i<\infty$ and $\sum_{k=0}^{\infty}\vert \psi(x^{k+1})-\psi(\tilde{x}^k)\vert<\infty$.
\label{truncation:finite move}
\end{truncation:finite move}

\begin{proof}
The estimate $\sum_{k=0}^{\infty}\Vert x^{k+1}-\tilde{x}^k\Vert\leq m\kappa\sum_{i=0}^{\infty}\epsilon_i<\infty$ follows directly from Definition \ref{def: truncation} and the settings in the truncation step. Hence, due to the local Lipschitz continuity of $\psi$ and the boundedness of $\{x^k\}_{k}$, there exists a constant $L_{\psi}>0$ such that $\sum_{k=0}^{\infty}\vert \psi(x^{k+1})-\psi(\tilde{x}^k)\vert\leq L_{\psi}m\kappa\sum_{i=0}^{\infty}\epsilon_i<\infty$.
\end{proof}

The next theorem is a weak global convergence result for Algorithm~\ref{alg:main}. In the proof, we combine our specific stepsize strategy and the truncation step to guarantee accuracy of the model and sufficient descent in $\psi$. More specifically, on the one hand, under the assumption that $\Vert g^k\Vert$ has a positive lower bound, the stepsize strategy ensures that $\Delta_k$ can not be arbitrarily small; on the other hand, the truncation step guarantees that the stepsize safeguard has a positive lower bound on an infinite set of iterations. Combining these two results, we would conclude that the total descent is infinite, which is a contradiction. We want to point out here that even if Algorithm~\ref{alg:main} may not need to compute $\rho^k_2$ in some steps, we still use it in our analysis.

\newtheorem{weak global convergence}[lem: partial psi]{Theorem}
\begin{weak global convergence}\label{liminf conv}
Suppose that the conditions (A.1) and (B.1)-(B.4) are satisfied and that $\psi$ can be truncated. Furthermore, let us assume that Algorithm~\ref{alg:main} does not terminate in finitely many steps and let $\{x^k\}_{k}$ be the generated sequence of iterates. Then, it holds that
\begin{equation*}
    \liminf_{k\rightarrow\infty}\Vert g^k\Vert=0.
\end{equation*}
\label{weak global convergence}
\end{weak global convergence}

\begin{proof}
Since Algorithm~\ref{alg:main} does not terminate in finitely many steps, we have
\begin{equation*}
    \psi'(x^k;d(x^k))<0 \quad \text{and} \quad \Vert d(x^k)\Vert=1 \quad \forall~k\in \mathbb{N}.
\end{equation*}
Suppose there exist $\epsilon>0$ and $K\in\mathbb{N}_+$ such that $\Vert g^k\Vert>\epsilon$ for all $k\geq K$. By the definition of the Cauchy point, we know that $m_k(s^k_C)\leq m_{k}(\Vert s^k\Vert d(x^k))$. Using \eqref{model reduce}, we obtain
\begin{equation*}
\begin{split}
    m_k(0)-m_k(s^k)\geq (1-\ell(\Vert s^k\Vert))[ m_k(0)-m_k(s^k_C)] \geq (1-\ell(\Vert s^k\Vert))[ m_k(0)-m_{k}(\Vert s^k\Vert d(x^k))],
\end{split}
\end{equation*}
and we have
\begin{equation*}
\begin{split}
\langle g^k, s^k-(1-\ell(\Vert s^k\Vert))\Vert s^k\Vert d(x^k)\rangle\leq -\frac{1}{2}\langle s^k,B^k  s^k\rangle+\frac{1-\ell(\Vert s^k\Vert)}{2}\langle \Vert s^k\Vert d(x^k),B^k \Vert s^k\Vert d(x^k)\rangle\leq \kappa_B\Vert s^k\Vert^2.
\end{split}
\end{equation*}
Due to $\bar{g}^k=-d(x^k)$, it holds that
\begin{equation*}
\langle -d(x^k),\bar{s}^k-(1-\ell(\Vert s^k\Vert))d(x^k)\rangle \leq \frac{\kappa_B}{\Vert g^k\Vert}\Vert s^k\Vert\leq \frac{\kappa_B}{\epsilon}\Vert s^k\Vert
\end{equation*}
for every $k\geq K$ which implies $-\langle d(x^k),\bar{s}^k\rangle<-1+\ell(\Vert s^k\Vert)+\frac{\kappa_B}{\epsilon}\Vert s^k\Vert$. Hence, we have
\begin{equation}
\Vert \bar{s}^k-d(x^k)\Vert=\sqrt{2-2\langle d(x^k),\bar{s}^k\rangle}\leq\sqrt{2 \ell(\Vert s^k\Vert)+\frac{2\kappa_B}{\epsilon}\Vert s^k\Vert}.\label{bound:s-d}
\end{equation}
Using the boundedness of $\{x^k\}_{k}$, the fact that $\psi'(x;d)$ is Lipschitz in $d$ for local $x$ (see, e.g., \citep{Nonsmooth-anal-opti-Clason-18}), and \eqref{bound:s-d}, we derive that there exists $L_\psi>0$ such that
\begin{equation*}
\vert \psi'(x^{k};\bar{s}^{k})-\psi'(x^k;d(x^k))\vert \leq L_\psi\sqrt{2 \ell(\Vert s^k\Vert)+\frac{2\kappa_B}{\epsilon}\Vert s^k\Vert} \quad \forall~k \geq K,
\end{equation*}
which combined with $g^k=u(x^k)d(x^k)$ and \eqref{u(x)} yields
\begin{equation}\label{bound:psi-g}
\begin{split}
    \alpha_k\psi'(x^k;\bar{s}^k)-\langle g^k,\alpha_k d(x^k)\rangle = & \alpha_k (\psi'(x^k;\bar{s}^k)-u(x^k))\\
    \leq & \alpha_k (\psi'(x^k;d(x^k))-u(x^k))+\alpha_k L_\psi\sqrt{2 \ell(\Vert s^k\Vert)+\frac{2\kappa_B}{\epsilon}\Vert s^k\Vert}\\
    \leq & \alpha_k L_\psi\sqrt{2 \ell(\Vert s^k\Vert)+\frac{2\kappa_B}{\epsilon}\Vert s^k\Vert}.
\end{split}
\end{equation}
Together with \eqref{bound:s-d}, we obtain
\begin{equation} \label{bound:g-s-d}
\begin{split}
    \langle g^k, d(x^k)\rangle-\langle g^k,\bar{s}^k\rangle  \leq \Vert g^k\Vert \Vert d(x^k)-\bar{s}^k \Vert \leq  \Vert g^k\Vert\sqrt{2 \ell(\Vert s^k\Vert)+\frac{2\kappa_B}{\epsilon}\Vert s^k\Vert}.
\end{split} 
\end{equation}
Using the definition of the function $h$ in \eqref{eq-B.3:global} and combining \eqref{bound:psi-g} and \eqref{bound:g-s-d}, it follows
\begin{align*}
   \psi(x^k+\alpha_k\bar{s}^k)-\psi(x^k)-\langle g^k,\alpha_k\bar{s}^k\rangle\leq & \alpha_k\left[(L_{\psi}+\Vert g^k\Vert)\sqrt{2 \ell(\Vert s^k\Vert)+\frac{2\kappa_B}{\epsilon}\Vert s^k\Vert}+h(\Delta_k)\right]\\
   \leq & \alpha_k\left[(L_{\psi}+\Vert g^k\Vert)\sqrt{2 \ell(\Delta_k)+\frac{2\kappa_B}{\epsilon}\Delta_k}+h(\Delta_k)\right]. 
\end{align*}
For all $k\geq K$, setting $\nu(\Delta_k) := [2\ell(\Delta_k)+ 2\kappa_B \Delta_k / \epsilon]^\frac12$, we now get
\begin{align*}
m_k(0)-m_k(\alpha_k\bar{s}^k) & =  \alpha_k\left[\langle -g^k,\bar{s}^k\rangle-\frac{\alpha_k}{2}\langle\bar{s}^k,B^k\bar{s}^k\rangle\right] 
\\ & \geq \alpha_k\left[\langle -g^k,d(x^k)\rangle-\Vert g^k\Vert\sqrt{2 \ell(\Vert s^k\Vert)+\frac{2\kappa_B}{\epsilon}\Vert s^k\Vert} -\frac{\alpha_k\kappa_B}{2}\right] \\ & \geq  \alpha_k\Vert g^k\Vert \left[1- \nu(\Delta_k) -\frac{\kappa_B}{2\epsilon}\Delta_k\right],
\end{align*}
and
\begin{align*}
1-\rho^2_k & =1-\frac{\psi(x^k)-\psi(x^k+\alpha_k\bar{s}^k)}{m_k(0)-m_k(\alpha_k\bar{s}^k)} =\frac{\psi(x^k+\alpha_k\bar{s}^k)-\psi(x^k)-\alpha_k \langle g^k,\bar{s}^k\rangle-\frac{\alpha_k^2}{2}\langle \bar{s}^k, B^k \bar{s}^k\rangle}{m_k(0)-m_k(\alpha_k\bar{s}^ks)}\\
& \leq \frac{(L_{\psi}+\Vert g^k\Vert)\nu(\Delta_k)+h(\Delta_k)+\frac{\kappa_B}{2}\Delta_k}{\Vert g^k\Vert \left[1-\nu(\Delta_k) -\frac{\kappa_B}{2\epsilon}\Delta_k\right]}  \leq\frac{(\frac{L_{\psi}}{\epsilon}+1)\nu(\Delta_k)+\frac{1}{\epsilon}h(\Delta_k)+\frac{\kappa_B}{2\epsilon}\Delta_k}{1-\nu(\Delta_k) -\frac{\kappa_B}{2\epsilon}\Delta_k}.    
\end{align*}
Thus, there exists $\sigma \in (0,{\epsilon}/{\kappa_B})$, such that for every $k\geq K$ with $\Delta_k\leq {\sigma}$ it holds that $1-\rho^2_k<1-\eta_1$. This implies $\rho^2_k>\eta_1$ for all $k\geq K$ satisfying $\Delta_k<{\sigma}$ which means that those steps are at least ``subsuccessful''. Hence, we can infer
\begin{equation}
\Delta_k\geq\min\{\Delta_K,r_1{\sigma}\} \quad \forall~k\geq K.\label{low_bound:Delta}
\end{equation}
Next, let us set $\mathcal{K}=\{k\geq K\mid\rho^1_k\geq\eta_1\ \text{or}\ \rho^2_k\geq\eta\}$ and 
\begin{equation*}
    \mathcal{K}_i=\{k\in\mathcal{K}\mid x^k\in N_i=S_{i}\backslash S_{i+1} \} \quad i=0,1,2,\cdots,m.
\end{equation*}
Due to \eqref{low_bound:Delta}, we have $|\mathcal{K}|=\infty$ and applying Lemma \ref{truncation:finite move}, it follows
\begin{equation*}
\begin{split}
    \sum_{k\in\mathcal{K}}[\psi(x^k)-\psi(\tilde{x}^k)] \leq & \sum_{k\in\mathcal{K}}[\psi(x^k)-\psi(x^{k+1})]+\sum_{k\in\mathcal{K}}\vert \psi(x^{k+1})-\psi(\tilde{x}^k)\vert\\
    \leq & \psi(x^0)-L_b+\sum_{k=0}^\infty\vert \psi(x^{k+1})-\psi(\tilde{x}^k)\vert< \infty,
\end{split}
\end{equation*}
where we used $x^{k+1}=\tilde{x}^k=x^k$ for all $k\notin\mathcal{K}$. Hence, we have
\begin{equation}
    \sum_{k\in\mathcal{K}}[m_k(0)-m_k(\tilde{x}^k-x^k)]\leq\frac{1}{\eta}\sum_{k\in\mathcal{K}}[\psi(x^k)-\psi(\tilde{x}^k)]<\infty.
\end{equation}
We now define the index $i_0=\max\{i \in \{0,1,2,\cdots,m\} : |\mathcal{K}_i|=\infty\} $. By the optimality of $i_0$, we can conclude that only finitely many elements of the sequence $\{x^k\}_{k}$ belong to $S_{i_0+1}$. This implies that the truncation operator $T$ is only applied a finite number of times on points in $N_{i_0}=S_{i_0}\backslash S_{i_0+1}$. In particular, $T$ will move points from $S_{i_0}\backslash S_{i_0+1}$ to the set $S_{i_0+1}$ and after a certain number of iterations $K'$, the counter $c_{i_0}$ will not be updated anymore, i.e., we have $c_{i_0} \equiv c$ for some $c$. Then, it follows $\Gamma(x^k)\geq \epsilon_{{c}}$ for all $x^k \in N_{i_0}$ and $k\geq K'$ and by (B.4) there exists $\epsilon'>0$ such that we have $\Gamma(x^k,\bar{s}^k)\geq\epsilon'$ for all $x^k\in N_{i_0}$ and $k\geq K'$. Combining \eqref{Cau_Reduce_step}, \eqref{model reduce0}, and \eqref{model_convex}, we can always guarantee descent in the model. In particular, we obtain
\begin{equation}
m_k(0)-m_k(\tilde{x}^k-x^k)\geq \frac{\Vert \tilde{x}^k-x^k\Vert}{4\Vert s^k\Vert}\cdot {\delta}_1\Vert g^k\Vert\min\left\{\Delta_k, {\delta}_2\Vert g^k\Vert\right\},\label{model_descent}
\end{equation}
where ${\delta}_1=\min\{\gamma_1,1\}$ and ${\delta_2}=\min\{\gamma_2,1/\kappa_B\}$.
Thus, we can conclude that
\begin{equation*}
\begin{split}
\infty&>\sum_{k\in\mathcal{K}}[m_k(0)-m_k(\tilde{x}^k-x^k)] 
\geq \sum_{k\in\mathcal{K}_{i_0},k\geq K'}\frac{{\delta}_1\Vert \tilde{x}^k-x^k\Vert}{4\Vert s^k\Vert}\Vert g^k\Vert\min\left\{\Delta_k, {\delta}_2\Vert g^k\Vert\right\} \\
&\geq\sum_{k\in\mathcal{K}_{i_0},k\geq K'}\frac{{\delta}_1\epsilon}{4}\frac{\min\left\{\Gamma(x^k,\bar{s}^k),\Vert s^k\Vert\right\}}{\Vert s^k\Vert}\min\left\{\Delta_K, r_1\sigma,{\delta}_2\epsilon\right\}\\
&\geq\sum_{k\in\mathcal{K}_{i_0},k\geq K'}\frac{{\delta}_1\epsilon}{4}\min\left\{\frac{\epsilon'}{{\Delta}_{\max}},1\right\}\min\left\{\Delta_K, r_1{\sigma},{\delta}_2\epsilon\right\}=\infty,    
\end{split}
\end{equation*}
which is a contradiction.
\end{proof}

 Next, we prove a stronger version of our global result under the additional assumption (A.2). Specifically, we show that every accumulation point of Algorithm~\ref{alg:main} is a stationary point of \eqref{obj fun}. This is a standard global convergence result, see, e.g., \cite{Nonsmooth-TR-Qi-Sun-94}. 

\newtheorem{strong global convergence}[lem: partial psi]{Theorem}
\begin{strong global convergence}
Let the conditions (A.1)--(A.2) and (B.1)--(B.4) be satisfied and suppose that $\psi$ can be truncated. Assume that Algorithm~\ref{alg:main} does not terminate after finitely many steps and that it generates a sequence $\{x^k\}_{k}$ with an accumulation point $x^*$. Then, $x^*$ is a stationary point of \eqref{obj fun}.\label{strong global convergence}
\end{strong global convergence}

\begin{proof}
We assume that $x^*$ is not a stationary point of \eqref{obj fun}. By (A.2) there exist $r,\epsilon>0$ such that $\Vert g(y)\Vert\geq \epsilon$ for all $y\in B_r(x^*)$. Let us set $A^k=\max\{\psi(x^{k+1})-\psi(\tilde{x}^k),0\}$ for $k\in\mathbb{N}$. Applying Lemma~\ref{truncation:finite move}, we know that $\sum_{k=0}^{\infty}A^k<\infty$. For any $k'>k$, we have
\begin{equation*}
\begin{split}
    \psi(x^{k'})& = \psi(x^k)+\sum_{t=k}^{k'-1}(\psi(x^{t+1})-\psi(\tilde{x}^t))+\sum_{t=k}^{k'-1}(\psi(\tilde{x}^t)-\psi(x^t))\\
    &\leq\psi(x^k)+\sum_{t=k}^{k'-1}A^t\leq \psi(x^k)+\sum_{t=k}^{\infty}A^t,
\end{split}
\end{equation*}
where we used the descent property $\psi(\tilde{x}^t)-\psi(x^t)\leq 0$. Consequently, we can infer
\begin{equation*}
    \limsup_{k'\rightarrow\infty}\psi(x^{k'})\leq \liminf_{k\rightarrow\infty}\psi(x^k)+\lim_{k\rightarrow\infty}\sum_{t=k}^{\infty}A^t\leq \liminf_{k\rightarrow\infty}\psi(x^k),
\end{equation*}
which implies that $\{\psi(x^k)\}_{k}$ converges. Next, Lemma~\ref{truncation:finite move} implies that there exists a constant $K\in\mathbb{N}$ such that $\sum_{k=K}^{\infty}\Vert x^{k+1}-\tilde{x}^k\Vert\leq \frac{r}{4}$. There is a subsequence $\{x^k\}_{k\in\mathcal{K}}\subseteq\{x^k\}_{k=K}^{\infty}$ satisfying $\{x^k\}_{k\in\mathcal{K}}\subseteq B_{r/4}(x^*)$ and $x^k\rightarrow x^*,\ k\rightarrow\infty, k\in\mathcal{K}$. For any $k\in\mathcal{K}$, since $\Vert g(y)\Vert\geq \epsilon,\ \forall~y\in B_r(x^*)$ and $\liminf_{k'} \Vert g^{k'}\Vert=0$ by Theorem~\ref{weak global convergence}, there must be some $k'\geq k$ such that $x^{k'}\notin B_{r}(x^*)$. Set $l(k)=\sup\{k'\geq k:x^t\in B_r(x^*),\ \forall~k\leq t\leq k'\}$. Thus, it holds that
\begin{equation*}
\begin{split}
    \sum_{t=k}^{l(k)}\Vert \tilde{x}^t-x^t\Vert+\frac{r}{4}&\geq\sum_{t=k}^{l(k)}(\Vert x^{t+1}-\tilde{x}^t\Vert+\Vert \tilde{x}^t-x^t\Vert)\\
    &\geq\sum_{t=k}^{l(k)}\Vert x^{t+1}-x^t\Vert\geq \Vert x^{l(k)+1}-x^k\Vert\geq\frac{3r}{4}
\end{split}
\end{equation*}
and it follows $\sum_{t=k}^{l(k)}\Vert \tilde{x}^t-x^t\Vert\geq \frac{r}{2}$. Mimicking the last steps in the proof of Theorem~\ref{weak global convergence}, we get
\begin{align*}
    \psi(x^k)-\psi(x^{l(k)+1})\geq& \sum_{t=k}^{l(k)}(\psi(x^t)-\psi(\tilde{x}^t)-A^t)\\
    \geq& \sum_{k\leq t\leq l(k),\ \rho^1_t\geq\eta_1 \text{ or } \rho^2_t\geq \eta} \eta[m_t(0)-m_t(\tilde{x}^t-x^t)]-\sum_{t=k}^{\infty }A^t\\
    \geq& \sum_{k\leq t\leq l(k),\ \rho^1_t\geq\eta_1 \text{ or } \rho^2_t\geq \eta} \frac{\eta{\delta}_1\Vert \tilde{x}^t-x^t\Vert}{4\Vert s^t\Vert}\Vert g^t\Vert\min\left\{\Delta_t, {\delta}_2\Vert g^t\Vert\right\}-\sum_{t=k}^{\infty }A^t\\
    \geq& \sum_{k\leq t\leq l(k),\ \rho^1_t\geq\eta_1 \text{ or } \rho^2_t\geq\eta} \frac{\eta{\delta}_1\epsilon\Vert \tilde{x}^t-x^t\Vert}{4}\min\left\{1,\frac{{\delta}_2\epsilon}{{\Delta}_{\max}}\right\}-\sum_{t=k}^{\infty }A^t\\
    \geq&  \frac{\eta{\delta}_1\epsilon r}{8}\min\left\{1, \frac{{\delta}_2\epsilon}{ {\Delta}_{\max}} \right\}-\sum_{t=k}^{\infty }A^t.
\end{align*}
Taking the limit $\mathcal{K}\ni k\rightarrow\infty$ we obtain the contradiction $0\geq \frac{\eta{\delta}_1\epsilon r}{8}\min\{1, {\delta}_2\epsilon{\Delta}_{\max}^{-1}\}$.
\end{proof}

\newtheorem{rmk-bdd-x}[lem: partial psi]{Remark}
\begin{rmk-bdd-x}
Theorem~\ref{strong global convergence} essentially establishes a similar result as in \cite[Theorem 3.4]{Nonsmooth-TR-Qi-Sun-94}. We notice that instead of boundedness of the level set $\{x\in\mathbb{R}^n\mid \psi(x)\leq \psi(x^0)\}$ (which was used in \cite{Nonsmooth-TR-Qi-Sun-94}), we need to work with the slightly stronger assumption (B.1) here since the truncation step can increase the objective function value $\psi$. However, if the function $\psi$ satisfies a Lipschitz-type assumption, i.e., if there are $\epsilon > 0$ and $L \geq 0$ such that $|\psi(x)-\psi(y)| \leq L\|x-y\|$ for all $x,y$ with $\|x-y\| \leq \epsilon$, then the proof of Lemma~\ref{truncation:finite move} implies $\sum_{k=0}^{\infty}\vert \psi(x^{k+1})-\psi(\tilde{x}^k)\vert<\infty$. This can be combined with the descent property $\psi(\tilde{x}^k)\leq \psi(x^k)$ of the trust-region step to show that the iterates $\{x^k\}$ will stay in the level set $\{x\in\mathbb{R}^n:\psi(x)\leq \zeta\}$ where $\zeta=\psi(x^0)+\sum_{k=0}^{\infty}\vert \psi(x^{k+1})-\psi(\tilde{x}^k)\vert<\infty$. Therefore, (B.1) can be substituted by a more classical level set condition in such a situation.
\label{rmk-bdd-x}
\end{rmk-bdd-x}

Finally, via utilizing the natural residual, it is possible to obtain strong lim-type convergence of Algorithm~\ref{alg:main}. Actually, compared with Theorem~\ref{liminf conv} which states that along a subsequence, $g^k$ converges to zero, the next theorem proves that some nonsmooth residual of $x^k$ converges to zero along the whole sequence.

\newtheorem{lim global convergence}[lem: partial psi]{Theorem}
\begin{lim global convergence}\label{lim global convergence}
Suppose that the same assumptions stated in Theorem~\ref{strong global convergence} are satisfied. Then, it holds that $\lim_{k\rightarrow\infty}\Vert F_{\mathrm{nat}}^{\Lambda}(x^k)\Vert=0$.
\end{lim global convergence}

\begin{proof}
Suppose that there exists $\epsilon>0$ and an infinite subsequence $\{x^k\}_{k\in\mathcal{K}}$ of $\{x^k\}_{k=0}^{\infty}$ satisfying
\begin{align}
   \Vert F_{\mathrm{nat}}^{\Lambda}(x^k)\Vert\geq\epsilon \quad  \forall~k\in\mathcal{K}. \label{lim glo_conv: contradict}
\end{align}
By (B.1),  $\{x^k\}_{k\in\mathcal{K}}$ has another subsequence $\{x^k\}_{k\in\mathcal{K}_1}$ with limit $x^*=\lim_{\mathcal{K}_1 \ni k\rightarrow\infty}x^k$. By Theorem~\ref{strong global convergence}, $x^*$ is a stationary point of \eqref{obj fun} with $F_{\mathrm{nat}}^{\Lambda}(x^*)=0$. Using the continuity of $F_{\mathrm{nat}}^{\Lambda}$, this contradicts \eqref{lim glo_conv: contradict}.
\end{proof}

\section{Fast Local Convergence}
\label{sec: loc_conv}
To the best of our knowledge, there are very limited local convergence results for nonsmooth trust-region type methods and most of the existing work only focuses on the global convergence analysis, see, e.g., \cite[][]{Nonsmooth-TR-Qi-Sun-94,Nonsmooth-TR-Riemannian-Grohs-16, TR-quadratic-model-Akbari-15}. In this section, we investigate local properties of our algorithm. Specifically, we will establish fast local convergence for the composite program $\psi=f+\varphi$ when $f$ is a smooth mapping and $\varphi$ is real-valued convex and partly smooth relative to an affine subspace. Our local results require that the first- and second-order information, i.e., $g^k$ and $B^k$, are chosen as the Riemannian gradient and the Riemannian Hessian with respect to some active manifold. We will also show that such information can be derived without knowing the active manifold under some suitable assumptions.

\subsection{Definitions and Assumptions}
In this subsection, we state some elementary definitions and assumptions. The family of partly smooth functions was originally introduced in \citep{Lewis-02} and plays a fundamental role in nonsmooth optimization. In particular, the concept of partly smoothness is utilized in the convergence analysis of nonsmooth optimization algorithms and to derive activity identification properties, see, e.g., \cite{partly-smooth-FBE-Liang-17, partly-smooth-SAGA-Poon-18}. Since in \eqref{obj fun:composite}, the mapping $\varphi$ is real-valued convex, we use the definition of partly smooth functions given in \citep{partly-smooth-FBE-Liang-17}. For a more general version and further details, we refer to \citep{Lewis-02}.

\newtheorem{def:part-smooth}[lem: partial psi]{Definition}
\begin{def:part-smooth}\cite[Definition 3.1]{partly-smooth-FBE-Liang-17}
A proper convex and lower semicontinuous function $\varphi$ is said to be partly smooth at $x$ relative to a set $\mathcal{M}$ containing $x$ if $\partial\varphi(x)\neq \emptyset$ and we have:
\begin{itemize}
    \item[{(i)}] \textbf{Smoothness}: $\mathcal{M}$ is a $C^2$-manifold around $x$ and $\varphi$ restricted to $\mathcal{M}$ is $C^2$ around $x$;
    \item[{(ii)}] \textbf{Sharpness}: The tangent space $T_{\mathcal{M}}(x)$ coincides with $T_x:=\text{par}(\partial\varphi(x))^{\perp}$, where $\text{par}(A)=\text{span}(A-A)$ for a convex set $A\subset \mathbb{R}^n$.
    \item[{(iii)}] \textbf{Continuity}: The set-valued mapping $\partial\varphi$ is continuous at $x$ relative to $\mathcal{M}$.
\end{itemize}\label{def:part-smooth}
\end{def:part-smooth}

Let $\{S_i\}_{i=0}^m$ be the sequence of sets associated with the truncation operator of $\psi$ and let $\{x^k\}_{k}$ be generated by Algorithm~\ref{alg:main}. We further consider an accumulation point $x^*$ of $\{x^k\}_k$ with $x^*\in S_{i^*}\backslash S_{i^*+1}$ and $i^*\in\{0,1,\cdots,m\}$ and we make the following assumptions.

\newtheorem{assump-part-smooth}[lem: partial psi]{Assumption}
\begin{assump-part-smooth}
 We consider the following conditions:
\begin{itemize}[leftmargin=7ex]
    \item[{(C.1)}\label{C.1}] The mapping $\varphi$ is partly smooth at $x^*$ relative to an affine subspace $\mathcal{M}$ and it holds that $B_r(x^*)\cap S_{i^*}=B_r(x^*)\cap\mathcal{M}$ for all $r\in(0,\Gamma(x^*))$. 
    \item[{(C.2)}\label{C.2}] The Riemannian Hessian $\nabla^2_{\mathcal{M}} \psi(x)$ is locally Lipschitz continuous around $x^*$ restricted to $\mathcal{M}$ and the second-order sufficient condition is satisfied at $x^*$, i.e., we have $\nabla^2_{\mathcal{M}}\psi(x^*)[\xi,\xi]\geq c\Vert \xi\Vert^2$ for some positive constant $c$ and all $\xi\in T_{\mathcal{M}}(x^*)$.
    \item[{(C.3)}\label{C.3}] The strict complementary condition $-\nabla f(x^*)\in\text{ri}\ \partial\varphi(x^*)$ is satisfied.
    \item[{(C.4)}\label{C.4}] For all $x\in S_i$, $y\in S_j$ with $i<j$, it holds that $\Gamma_{\max}\left(x,\overline{y-x}\right)\leq \Vert y-x\Vert$ where $\overline{y-x}=\frac{y-x}{\Vert y-x\Vert}$. For every $r\in(0,\Gamma(x^*))$, there exists $\epsilon(r) > 0$ such that $\Gamma(x)\geq \epsilon(r)$ for all  $x\in B_r(x^*)\cap\mathcal{M}$.
    \item[{(C.5)}\label{C.5}] The sequence $\{x^k\}_{k}$ converges with limit $\lim_{k\rightarrow\infty}x_k=x^*$.
\end{itemize}
\label{assump-part-smooth}
\end{assump-part-smooth}

Besides the partly smoothness, assumption (C.1) states that the local structure of $S_{i^*}$ around $x^*$ has to be affine.
The conditions (C.2), (C.3), and (C.5) are standard assumptions for finite active identification and have been used to establish local convergence rates. For instance, they appeared in \cite[][]{partly-smooth-FBE-Liang-17,partly-smooth-SAGA-Poon-18}.

In order to illustrate assumption (C.4), we consider the example $\varphi(x)=\Vert x\Vert_1$. Suppose that $x\in S_i,\ y\in S_j$ are two given points with $i<j$. Since $y$ has more zero-components as $x$ there exists a point on the line connecting $x$ and $y$ at which $\varphi$ is not differentiable. This immediately leads to $\Gamma_{\max}\left(x,\overline{y-x}\right)\leq \Vert y-x\Vert$. The second part in (C.4) requires that  $\Gamma$ does not decay sharply around $x^*$ restricted to $\mathcal{M}$. We will use condition (C.4) in the analysis of the truncation step.

\subsection{Riemannian Gradient and Riemannian Hessian} We now choose $g^k=F_{\mathrm{nor}}^{\Lambda}(\tau(x^k))$. The next lemma shows that this choice actually coincides with the Riemannian gradient of $\psi$ when $x^k$ lies in the manifold $\mathcal{M}$ and is close to $x^*$.

\newtheorem{rie_grad}[lem: partial psi]{Lemma}
\begin{rie_grad}
Suppose that the conditions (C.1) and (C.3) hold and that $x^*$ is a stationary point. Let $x\in B_r(x^*)\cap\mathcal{M}$ be given for $r\in(0,\Gamma(x^*))$ sufficiently small. Then, we have $F_{\mathrm{nor}}^{\Lambda}(\tau(x))=\psi'(x;d_s(x))d_s(x)=\nabla_{\mathcal{M}} \psi(x)$, where $\nabla_{\mathcal{M}}\psi(x)$ denotes the Riemannian gradient of $\psi$.\label{rie_grad}
\end{rie_grad}

\begin{proof}
By the stationarity of $x^*$ and (C.3), it is easy to see that $\nabla_{\mathcal{M}}\psi(x^*)=0\in\nabla f(x^*)+\text{ri}\ \partial\varphi(x^*)=\text{ri}\ \partial\psi(x^*).$ Thus, applying \cite[Corollary 21]{Geo-Aris-06} for $x\in\mathcal{M}$ close to $x^*$ and Lemma \ref{lem: partial psi}, we can conclude that $\nabla_{\mathcal{M}}\psi(x)=\mathbf{P}_{\partial\psi(x)}(0)=\psi'(x;d_s(x))d_s(x)=F_{\mathrm{nor}}^{\Lambda}(\tau(x))$.
\end{proof}

Since $g^k$ coincides with the Riemannian gradient, we naturally would like to choose $B^k$ as the associated Riemannian Hessian of $\psi$. We now show that this Hessian can be derived without knowing the underlying manifold $\mathcal{M}$ if we additionally assume that $\varphi$ is polyhedral. Specifically, the following lemma establishes a connection between the derivative of $F_{\mathrm{nor}}^{\Lambda}(\tau(x))$ and the Riemannian Hessian.

\newtheorem{rie_hess}[lem: partial psi]{Lemma}
\begin{rie_hess}
Suppose that the assumptions stated in Lemma~\ref{rie_grad} hold and that $\varphi$ is a polyhedral function. For $r\in(0,\Gamma(x^*))$ sufficiently small and all $x\in B_r(x^*)\cap\mathcal{M}$, it follows 
\[ V\mathcal{D} F_{\mathrm{nor}}^{\Lambda}(z)=\nabla^2_{\mathcal{M}} \psi(x), \quad \Lambda=\lambda I, \; \lambda>0, \] 
where $\mathcal{D}$ is the differential operator, $z=\tau(x)$, $V=\mathcal{D}\mathrm{prox}_{\varphi}^{\Lambda}(z)$ is the derivative of $\mathrm{prox}_{\varphi}^{\Lambda}$, and $\nabla^2_{\mathcal{M}} \psi(x)$ is the Riemannian Hessian. \label{rie_hess}
\end{rie_hess}

\begin{proof}
For $x\in\mathcal{M}$ and near $x^*$, by Definition \ref{def:part-smooth} and \cite[Fact 3.3]{partly-smooth-FBE-Liang-17}, we can decompose $\partial\varphi(x)$ as $\partial\varphi(x)=\left\{\nabla_{\mathcal{M}}\varphi(x)\right\}+\partial_{\mathcal{M}}^{\perp}\varphi(x)$, where $\partial_{\mathcal{M}}^{\perp}\varphi(x)\subseteq  T_{\mathcal{M}}(x)^{\perp}$. We can see that both $\nabla_{\mathcal{M}}\varphi(x)$ and $\partial_{\mathcal{M}}^{\perp}\varphi(x)$ restricted to $\mathcal{M}$ are continuous around $x^*$. Moreover, we have the decomposition $\nabla f(x)=\nabla_{\mathcal{M}}f(x)+\nabla_{\mathcal{M}}^{\perp} f(x)$ where $\nabla_{\mathcal{M}}^{\perp} f(x)\in T_{\mathcal{M}}(x)^{\perp}$.

Condition (C.3) implies $-\nabla_{\mathcal{M}}^{\perp} f(x^*)\in \text{ri}\ \partial_{\mathcal{M}}^{\perp}\varphi(x^*)$, which combined with part (ii) in Definition \ref{def:part-smooth} and the continuity of $\nabla_{\mathcal{M}}^{\perp} f \big|_{\mathcal{M}}$ and $\partial_{\mathcal{M}}^{\perp}\varphi\big|_{\mathcal{M}}$ leads to $-\nabla_{\mathcal{M}}^{\perp} f(x)\in \text{ri}\ \partial_{\mathcal{M}}^{\perp}\varphi(x)$, i.e., $0\in \left\{\nabla_{\mathcal{M}}^{\perp} f(x)\right\}+\text{ri}\ \partial_{\mathcal{M}}^{\perp}\varphi(x)$ for all $x\in B_r(x^*)\cap\mathcal{M}$ where $r>0$ is sufficiently small. Thus, it follows
\begin{equation*}
\nabla f(x)+\Lambda(z-x)=F_{\mathrm{nor}}^{\Lambda}(z)=\nabla_{\mathcal{M}} f(x)+\nabla_{\mathcal{M}} \varphi(x)\in\nabla f(x)+\text{ri}\ \partial\varphi(x),
\end{equation*}
which implies $\Lambda(z-x)\in \text{ri}\ \partial\varphi(x)$. Since $\varphi$ is polyhedral, the subdifferential $\partial\varphi(x)$ is locally constant around $x^*$ on $\mathcal{M}$. For any $d_1\in T_{\mathcal{M}}$, $d_2\in T_{\mathcal{M}}^{\perp}$ with $\Vert d_1\Vert,\Vert d_2\Vert$ sufficiently small, it holds that $\Lambda((z+d_1+d_2)-(x+d_1))=\Lambda(z-x)+\Lambda d_2\in\partial\varphi(x)=\partial\varphi(x+d_1)$, which implies that $\mathrm{prox}_{\varphi}^{\Lambda}(z+d_1+d_2)=x+d_1=\mathrm{prox}_{\varphi}^{\Lambda}(z)+d_1$. Consequently, we have $\mathcal{D}\mathrm{prox}_{\varphi}^{\Lambda}(z) = V =\mathbf{P}$, where $\mathbf{P}$ is the orthogonal projection operator onto $ T_{\mathcal{M}}$. The derivative of the normal map is $\mathcal{D} F_{\mathrm{nor}}^{\Lambda}(z)=\nabla^2 f(x)\mathbf{P}+\Lambda(I-\mathbf{P})$, which combined with the local linearity of $\varphi$ yields $\mathbf{P}\mathcal{D} F_{\mathrm{nor}}^{\Lambda}(z)=\mathbf{P}\nabla^2 f(x)\mathbf{P}=\nabla^2_{\mathcal{M}} \psi(x)$.
\end{proof}

\subsection{Convergence Analysis} We have the following finite active identification result.

\newtheorem{finite-identify-part-smooth}[lem: partial psi]{Lemma}
\begin{finite-identify-part-smooth}
Suppose that the assumptions in Theorem~\ref{strong global convergence} are satisfied and that the conditions (C.4)--(C.5) hold. Then for every $r\in(0,\Gamma(x^*))$ there exist infinitely many $k\in\mathbb{N}$ with $x^k\in B_r(x^*)\cap S_{i^*}$. \label{finite-identify-part-smooth}
\end{finite-identify-part-smooth}

\begin{proof}
Without loss of generality, we can assume $\{x^k\}_{k} \subseteq B_r(x^*)$. Let us set
\begin{equation*}
    \mathcal{K}_i=\{k\in\mathbb{N}\mid x^k\in S_i\backslash S_{i+1}\},\quad i=0,1,\cdots,m.
\end{equation*}
By assumption (C.4), for every $y\in B_r(x^*)$ with $y\in S_{i^*+1}$ we have
\begin{equation*}
    \Vert y-x^*\Vert\geq\Gamma_{\max}\left(x^*,\overline{y-x^*}\right)\geq \Gamma(x^*)>r>\Vert y-x^*\Vert,
\end{equation*}
which is a contradiction. Hence, it follows $|\mathcal{K}_i|=0$ for every $i>i^*$, i.e., $B_r(x^*)\cap S_{i^*+1} = \emptyset$.

Set $i_0=\max\{ i=0,1,\cdots,i^*\mid|\mathcal{K}_i|=\infty \}$. If $i_0\leq i^*-1$, since truncations on points in $S_{i_0}\backslash S_{i_0+1}$ only happen in a finite number of times, $\{\Gamma(x^k)\mid k\in\mathcal{K}_{i_0}\}$ has a positive lower bound, i.e., $\beta:=\inf_{k\in\mathcal{K}_{i_0}} \Gamma(x^k)>0$. For $k\in\mathcal{K}_{i_0}$ and by (C.4), we can conclude that $\Gamma(x^k)\leq \Gamma_{\max}(x^k,\overline{x^*-x^k})\leq \Vert x^k-x^* \Vert$. Using $x^k\rightarrow x^*$, ($\mathcal{K}_{i_0} \ni k\rightarrow\infty$), we can infer
\begin{equation*}
    0<\beta\leq \liminf_{k\in\mathcal{K}_{i_0},\ k\rightarrow\infty}\Gamma(x^k)\leq\lim_{k\in\mathcal{K}_{i_0},\ k\rightarrow\infty}\Vert x^k-x^*\Vert=0,
\end{equation*}
which is a contradiction. Thus, we have $i_0=i^*$, which finishes the proof.
\end{proof}

At the end of this subsection, we establish the local convergence rate by connecting our algorithm with a Riemannian trust-region method.

\newtheorem{loc-conv-part-smooth}[lem: partial psi]{Theorem}
\begin{loc-conv-part-smooth}
Suppose that the assumptions in Theorem~\ref{strong global convergence} hold and that the conditions (C.1)--(C.5) are satisfied. Furthermore, if for some sufficiently small $r\in(0,\Gamma(x^*))$ and every $k$ with $x^k\in B_r(x^*)\cap \mathcal{M}$, we choose $g^k=\nabla_{\mathcal{M}}\psi(x^k)$, $B^k=\nabla^2_{\mathcal{M}}\psi(x^k)$, and solve the trust-region subproblem exactly with solution $s^k \in T_{\mathcal{M}}(x^k)$, then $\{x^k\}_{k}$ converges to $x^*$ q-quadratically.
\label{loc-conv-part-smooth}
\end{loc-conv-part-smooth}

\begin{proof}
We can assume $\{x^k\}_{k} \in B_r(x^*)\cap \mathcal{M}$ for some small enough $r>0$, which combined with (C.1) and the fact $B_r(x^*)\cap S_{i^*+1}=\emptyset$ implies that there is no truncation. Since the trust-region subproblem is solved exactly in $T_{\mathcal{M}}(x^k)$ and $\mathcal{M}$ is affine, condition (C.2) can be utilized to show that the first acceptance test is always locally successful and hence the algorithm always skips the second acceptance mechanism. A detailed proof of this observation, which is also applicable in our situation, can be found in \cite[Theorem 4.9]{Num-Opti-Nocedal-06}. We can then infer that our algorithm locally coincides with a Riemannian trust-region method or a classical trust-region method in the tangent space $T_{\mathcal{M}}(x^*)$ and the trust-region radius eventually becomes inactive. Thus, the local quadratic convergence rate is achieved by following \cite[Chapter 7]{absil-08} or \cite[Theorem 4.9]{Num-Opti-Nocedal-06}.
\end{proof}

\newtheorem{rmk-loc-conv}[lem: partial psi]{Remark}
\begin{rmk-loc-conv}
Lemmas~\ref{rie_grad} and \ref{rie_hess} guarantee that we can set $g^k=\nabla_{\mathcal{M}}\psi(x^k)$ and $B^k=\nabla^2_{\mathcal{M}}\psi(x^k)$ without explicitly knowing $\mathcal{M}$ and that there is a (globally optimal) solution of the trust-region subproblem located in $T_\mathcal{M}(x^k)$. This solution actually has the minimal $\ell_2$-norm among all solutions. Since $g^k$ and $B^k$ operate on the tangent space of the active manifold, some practical algorithms, such as the CG-Steihaug method, can indeed recover $s^k$ in $T_\mathcal{M}(x^k)$, which leads to $x^k+s^k\in \mathcal{M}$. Therefore, given $x^k\in B_r(x^*)\cap \mathcal{M}$ for $k$ sufficiently large, we would have $x^{k+1}\in B_r(x^*)\cap \mathcal{M}$.
\label{rmk-loc-conv}
\end{rmk-loc-conv}

\section{Preliminary Numerical Results}
\label{numerics}

In this section, we test the efficiency of our proposed nonsmooth trust-region method by applying it to convex and nonconvex $\ell_1$-minimization problems. All numerical experiments are performed in MATLAB R2020a on a laptop with Intel(R) Core(TM) i7-7700HQ CPU @ 2.80GHz and 16GB memory.

We apply our framework to the $\ell_1$-minimization problem 
\begin{equation}
\min_{x\in\mathbb{R}^n}~f(x)+\mu{\Vert x\Vert_1},\label{problem-num}
\end{equation}
where $f:\mathbb{R}^n\rightarrow\mathbb{R}$ is a smooth function. Setting $\varphi(x)=\mu{\Vert x\Vert_1}$ and using \eqref{g(x):2}, we choose
\begin{equation}
g(x)=\lambda F_{\mathrm{nat}}^{\lambda}(x):=\lambda [x-\mathrm{prox}_{\varphi}^{\lambda}(x-\lambda^{-1}\nabla f(x))], \quad \lambda > 0.
\end{equation}
where the proximity operator is given explicitly by $(\mathrm{prox}_{\varphi}^{\lambda}(x))_i=\text{sign}(x_i)\max(\vert x_i\vert-\lambda^{-1}\mu,0)$. We now construct an element $M(x)\in\partial \mathrm{prox}_{\varphi}^{\lambda}(x-\lambda^{-1}\nabla f(x))$ as follows: $M(x)\in\mathbb{R}^{n\times n}$ is a diagonal matrix with diagonal entries
\begin{equation*}
(M(x))_{ii}=\begin{cases}1,&\text{if }\vert (x-\lambda^{-1}\nabla f(x))_i\vert >\lambda^{-1}\mu,\\
0,&\text{otherwise}.\end{cases}
\end{equation*}
Thus, $J(x)=I-M(x)(I-\lambda^{-1}\nabla^2 f(x))$ is a possible generalized Jacobian of $F_{\mathrm{nat}}^{\lambda}(x)$. Let us define the index sets
\begin{equation*}
\mathcal{I}(x):=\{i\in\{1,2,\dots,n\}: \vert (x-\lambda^{-1}\nabla f(x))_i\vert >\lambda^{-1}\mu\},
\end{equation*}
and
\begin{equation*}
    \mathcal{O}(x):=\{i\in\{1,2,\dots,n\}: \vert (x-\lambda^{-1}\nabla f(x))_i\vert \leq\lambda^{-1}\mu\}.
\end{equation*}
Then, $J(x)$ can be written in an alternative format:
\begin{equation}
J(x)=\begin{pmatrix}
\lambda^{-1}(\nabla^2 f(x))_{\mathcal{I}(x)\mathcal{I}(x)} & \lambda^{-1}(\nabla^2 f(x))_{\mathcal{I}(x)\mathcal{O}(x)}\\
0& I
\end{pmatrix}.
\end{equation}
In the following, we choose $B^k=\lambda J(x^k)$. For simplicity, we do not check the condition \eqref{model_convex}.

\subsection{The Lasso Problem}
\label{sec:lasso}

We first consider the Lasso problem where $f$ is a convex quadratic function
\begin{equation*}
f(x)=\frac{1}{2}\Vert Ax-b\Vert^2,    
\end{equation*}
and $b\in\mathbb{R}^m$ and $A=\mathbb{R}^{m\times n}$ are given.

It can be shown that $J(x)$ is positive semidefinite if $\lambda$ is sufficiently large \citep{Semismooth-composite-Xiao-18}. When solving the trust-region subproblem \eqref{sub-problem} and similar to the method presented in Appendix~\ref{solution_sub}, we first choose a suitable regularization parameter $t_k\geq 0$ and solve the linear system
\begin{equation}
(J^k+t_k I)p^k=- F_{\mathrm{nat}}^{\lambda}(x^k), \quad J^k=J(x^k), \label{Newton-sys}
\end{equation}
and then project $p^k$ onto the trust region, i.e., $s^k=\min\{\Delta_k,\Vert p^k\Vert\}{\bar p^k}$. Setting $g^k = g(x^k)$, $\mathcal{I}^k=\mathcal{I}(x^k)$, and $\mathcal{O}^k=\mathcal{O}(x^k)$, the linear system \eqref{Newton-sys} is equivalent to
\begin{equation*}
(1+t_k)p^k_{\mathcal{O}^k} =-g^k_{\mathcal{O}^k}, \quad 
(\lambda^{-1} (A^T A)_{\mathcal{I}^k\mathcal{I}^k}+t^k I)p^k_{\mathcal{I}^k}+\lambda^{-1} (A^T A)_{\mathcal{I}^k\mathcal{O}^k}p^k_{\mathcal{O}^k}  =-g^k_{\mathcal{I}^k},
\end{equation*}
which leads to
\begin{equation*}
p^k_{\mathcal{O}^k}=-\frac{1}{(1+t_k)}g^k_{\mathcal{O}^k}, \quad 
(\lambda^{-1} (A^T A)_{\mathcal{I}^k\mathcal{I}^k}+t^k I)p^k_{\mathcal{I}^k}=-g^k_{\mathcal{I}^k}-\lambda^{-1} (A^T A)_{\mathcal{I}^k\mathcal{O}^k}p^k_{\mathcal{O}^k}.
\end{equation*}
The second system is symmetric and can be much smaller than the original problem \eqref{Newton-sys}. It can be solved efficiently by applying the CG method.

Our test framework follows \citep{BecBobCan11,Semismooth-L1-Milzarek-14}:
\begin{itemize}
\item A sparse solution $\hat{x}\in\mathbb{R}^n$ with $n=512^2=262144$ is generated randomly with $k=[n/40]$ zero entries. The nonzero components are chosen from $\{1,2,\cdots,n\}$ uniformly with values given by $\hat{x}_i=\eta_1(i) 10^{{d \eta_2(i)}/{20}}$. Here, $\eta_1(i)$ and $\eta_2(i)$ are distributed uniformly in $\{-1,1\}$ and $[0,1]$, respectively and $d$ is a dynamic range. 
\item We randomly choose $\mathcal{J}\subseteq\{1,2,\cdots,n\}$ with $\vert \mathcal{J}\vert=m = n/8=32768$. The linear operator $A:\mathbb{R}^n\rightarrow\mathbb{R}^m$ is then defined via $Ax=(\texttt{dct}(x))_{\mathcal{J}}$ where $\texttt{dct}$ denotes the discrete cosine transform.
\item We set $b=A\hat{x}+\epsilon$ where $\epsilon\in\mathbb{R}^m$ is Gaussian noise with covariance matrix $\hat{\sigma}I_{m\times m}$, $\hat{\sigma}=0.1$. 
\end{itemize}

Given a tolerance $\epsilon$, we terminate whenever the condition $\lambda\left\Vert F_{\mathrm{nat}}^{\lambda}(x)\right\Vert\leq \epsilon$ is satisfied. In the algorithm, $\lambda$ is chosen adaptively to estimate the local Lipschitz constant of $\nabla f$, i.e., we have $\lambda = \lambda_k = \max\{10^{-3},\min\{\|x^{k+1}-x^k\|/\|\nabla f(x^{k+1})-\nabla f(x^k)\|,10^{3}\}$ if the step was successful. We compare our nonsmooth trust-region method (NTR) with the adaptive semi-smooth Newton (ASSN) method in \citep{Semismooth-composite-Xiao-18} and the 
fast iterative shrinkage-thresholding algorithm (FISTA) \cite{FISTA} for different tolerances $\epsilon\in\{10^0,10^{-1},10^{-2},10^{-4},10^{-6}\}$ and dynamic ranges $d\in\{20,40,60,80\}$. We report the average CPU time (in seconds) as well as the average number of $A$- and $A^T$-calls $N_A$ over 10 independent trials. 

The numerical comparisons are shown in Tables~\ref{dyna20}-\ref{dyna80}. From those results we can see that the nonsmooth trust-region method outperforms the first-order method FISTA and is quite competitive with the second-order method ASSN. Even if the second acceptance test and the stepsize safeguard is required to guarantee theoretical convergence, in the numerical experiments we find that with a group of suitably chosen parameters, our algorithm rarely or never fails in the first acceptance test and hence employs the mechanism of the stepsize safeguard, which prevents additional costs. The similar behavior of NTR and ASSN may stem from the fact that we utilize similar strategies as in ASSN \cite{Semismooth-composite-Xiao-18} for certain parameters.
Our results on $N_A$ are comparable with ASSN's results and are sometimes better. Because each of our iterations involves potential acceptance tests and truncation steps, our method overall requires slightly more CPU time to converge than ASSN.

\begin{table}[t]\setlength{\tabcolsep}{4pt} 
\centering
\caption{Numerical results with dynamic range 20 dB}\label{dyna20}
\begin{tabular}{lllllllllll}
\hline
& \multicolumn{2}{c}{$\epsilon:10^0$}& \multicolumn{2}{c}{$\epsilon:10^{-1}$}& \multicolumn{2}{c}{$\epsilon:10^{-2}$}& \multicolumn{2}{c}{$\epsilon:10^{-4}$}& \multicolumn{2}{c}{$\epsilon:10^{-6}$}\\
\cline{2-11}
& time & $N_A$ & time & $N_A$ & time & $N_A$ & time & $N_A$ & time & $N_A$\\
\hline
NTR  & 0.8275 & 86.8 & 1.2226 & 132.8 & 1.5549 & 172 & 2.0399 & 227.4 & 2.4704 & 280.6 \\
\hline
ASSN  & 0.7368 & 89.8 & 1.1409 & 145& 1.3583 & 173 & 1.9094 & 246.4 & 2.2844 & 298.2 \\
\hline
FISTA & 0.5337 & 59 & 1.3959 & 153 & 3.2304 & 353.2 & 13.9021 & 1490.4 & 33.7451 & 3581.2 \\
\hline
\end{tabular}
\end{table}

\begin{table}[t]\setlength{\tabcolsep}{4pt} 
\centering
\caption{Numerical results with dynamic range 40 dB}\label{dyna40}
\begin{tabular}{lllllllllll}
\hline
& \multicolumn{2}{c}{$\epsilon:10^0$}& \multicolumn{2}{c}{$\epsilon:10^{-1}$}& \multicolumn{2}{c}{$\epsilon:10^{-2}$}& \multicolumn{2}{c}{$\epsilon:10^{-4}$}& \multicolumn{2}{c}{$\epsilon:10^{-6}$}\\
\cline{2-11}
& time & $N_A$ & time & $N_A$ & time & $N_A$ & time & $N_A$ & time & $N_A$\\
\hline
NTR  & 1.7293 & 176.4 & 2.6661 & 280.4 & 3.1330 & 330.8 & 3.7463 & 402.2 & 4.1825 & 464.2 \\
\hline
ASSN  & 1.5227 & 182.2 & 2.3414 & 285.4 & 2.7751 & 338.6 & 3.3216 & 407 & 3.6687 & 459.2 \\
\hline
FISTA  & 2.2007 & 234.8 & 3.9217 & 418.4 & 7.7210 & 817 & 26.6733 & 2804.6 & 57.0018 & 5991.6 \\
\hline
\end{tabular}
\end{table}

\begin{table}[t]\setlength{\tabcolsep}{4pt} 
\centering
\caption{Numerical results with dynamic range 60 dB}\label{dyna60}
\begin{tabular}{lllllllllll}
\hline
& \multicolumn{2}{c}{$\epsilon:10^0$}& \multicolumn{2}{c}{$\epsilon:10^{-1}$}& \multicolumn{2}{c}{$\epsilon:10^{-2}$}& \multicolumn{2}{c}{$\epsilon:10^{-4}$}& \multicolumn{2}{c}{$\epsilon:10^{-6}$}\\
\cline{2-11}
& time & $N_A$ & time & $N_A$ & time & $N_A$ & time & $N_A$ & time & $N_A$\\
\hline
NTR  & 2.9135 & 303.4 & 3.7794 & 398.4 & 4.4177 & 471.2 & 5.1089 & 562.4 & 5.6559 & 632 \\
\hline
ASSN  & 2.4508 & 295.4 & 3.4378 & 416.4 & 4.0051 & 492 & 4.6492 & 582.4 & 5.1121 & 642.4 \\
\hline
FISTA  & 5.9349 & 630.6 & 9.0296 & 951.6 & 14.7183 & 1548 & 39.5980 & 4164.2 & 79.6442 & 8355.4 \\
\hline
\end{tabular}
\end{table}

\begin{table}[t]\setlength{\tabcolsep}{4pt} 
\centering
\caption{Numerical results with dynamic range 80 dB}\label{dyna80}
\begin{tabular}{lllllllllll}
\hline
& \multicolumn{2}{c}{$\epsilon:10^0$}& \multicolumn{2}{c}{$\epsilon:10^{-1}$}& \multicolumn{2}{c}{$\epsilon:10^{-2}$}& \multicolumn{2}{c}{$\epsilon:10^{-4}$}& \multicolumn{2}{c}{$\epsilon:10^{-6}$}\\
\cline{2-11}
& time & $N_A$ & time & $N_A$ & time & $N_A$ & time & $N_A$ & time & $N_A$\\
\hline
NTR  & 3.6748 & 411 & 5.3121 & 614 & 5.9883 & 702.2 & 6.5791 & 803.6 & 7.2349 & 868.8 \\
\hline
ASSN  & 3.6290 & 482.8 & 4.5779 & 601 & 4.9879 & 690.6 & 5.7141 & 780.6 & 6.3010 & 865.4 \\
\hline
FISTA  & 20.9653 & 2222.4 & 25.6927 & 2673.2 & 33.2373 & 3527 & --- & --- & ---- & ---- \\
\hline
\end{tabular}
\end{table}

\subsection{Nonconvex Binary Classification}

We consider a second, nonconvex binary classification problem, \cite{mason1999boosting, wang2017stochastic}, where $f$ is given as follows:
\begin{equation*}
    f(x)=\frac{1}{N}\sum_{i=1}^N\left(1-\tanh(b_i\cdot a_i^T x)\right).
\end{equation*}
Here, the datapoints $a_i\in \mathbb{R}^n$ and features $b_i\in\{\pm 1\}$ are selected from the datasets \texttt{CINA} ($N = 16033$, $n = 132$) and \texttt{gisette} ($N = 6000$, $n = 5000$). We set $\mu=0.01$.

Although positive semidefiniteness of $J^k$ is not guaranteed in the nonconvex case, we reuse the method described in Section~\ref{sec:lasso} to solve the trust-region subproblem with $t_k=0$. The parameter $\lambda$ is updated adaptively as before.  We also use the same stopping criterion as in Section~\ref{sec:lasso}. We compare our nonsmooth trust-region method (NTR) with the stochastic semismooth Newton method with variance reduction (S4N-VR) \cite{Milzarek-19-sto-semismooth} for different tolerances $\epsilon\in\{10^0,10^{-1},10^{-2},10^{-4},10^{-6}\}$.

The CPU time (in seconds) and the number of $A$- and $A^T$-calls $N_A$ are reported in Table~\ref{CINA} and Table~\ref{gisette}, respectively. The results for S4N-VR are averaged over 10 independent trials while the results for NTR are based on one (deterministic) trial. While S4N-VR achieves slightly better results on \texttt{CINA}, NTR outperforms S4N-VR on the second dataset \texttt{gisette}.

Although the performance of NTR is still not perfect, our preliminary results underline that the proposed class of nonsmooth trust-region methods is promising and allows us to handle nonsmooth nonconvex optimization problems from a different perspective.

\begin{table}[t]\setlength{\tabcolsep}{4pt} 
\centering
\caption{Numerical results for \texttt{CINA}}\label{CINA}
\begin{tabular}{lllllllllll}
\hline
& \multicolumn{2}{c}{$\epsilon:10^0$}& \multicolumn{2}{c}{$\epsilon:10^{-1}$}& \multicolumn{2}{c}{$\epsilon:10^{-2}$}& \multicolumn{2}{c}{$\epsilon:10^{-4}$}& \multicolumn{2}{c}{$\epsilon:10^{-6}$}\\
\cline{2-11}
& time & $N_A$ & time & $N_A$ & time & $N_A$ & time & $N_A$ & time & $N_A$\\
\hline
NTR  & 0.03567 & 2.5995 & 0.08788 & 9.6904 & 0.2794 & 55.9480 & 0.2806 & 63.3571 & 0.3048 & 67.3268 \\
\hline
S4N-VR  & 0.02843 & 2.0983 & 0.08124 & 6.0051 & 0.2090 & 25.8886 & 0.2422 & 28.3287 & 0.3722 & 45.6214 \\
\hline
\end{tabular}
\end{table}

\begin{table}[t]\setlength{\tabcolsep}{4pt} 
\centering
\caption{Numerical results for \texttt{gisette}}\label{gisette}
\begin{tabular}{lllllllllll}
\hline
& \multicolumn{2}{c}{$\epsilon:10^0$}& \multicolumn{2}{c}{$\epsilon:10^{-1}$}& \multicolumn{2}{c}{$\epsilon:10^{-2}$}& \multicolumn{2}{c}{$\epsilon:10^{-4}$}& \multicolumn{2}{c}{$\epsilon:10^{-6}$}\\
\cline{2-11}
& time & $N_A$ & time & $N_A$ & time & $N_A$ & time & $N_A$ & time & $N_A$\\
\hline
NTR  & 0.3136 & 14.1908 & 0.5044 & 21.3416 & 0.6588 & 34.7592 & 0.9333 & 57.2328 & 1.0780 & 73.7344 \\
\hline
S4N-VR  & 0.9177 & 6.8913 & 1.6120 & 20.4010 & 2.3202 & 32.6835 & 5.2769 & 84.0567 & 6.7880 & 121.4267 \\
\hline
\end{tabular}
\end{table}

\section{Conclusion}
\label{sec: conclude}
 In this paper, we investigate a trust-region method for nonsmooth nonconvex optimization problems. In the proposed framework, the model functions are quadratic and cheap descent directions can be utilized. This allows us to construct cheap model functions and to apply standard algorithms for solving the resulting trust-region subproblems. We propose a novel combination of a stepsize safeguard for ensuring the accuracy of the model and an additional truncation step to enlarge the stepsize safeguard and to accelerate convergence. We present a detailed discussion of the global convergence properties under suitable and mild assumptions. In the case of composite-type problems, we also show that our method converges locally with a quadratic rate after the finite identification of the active manifold when the nonsmooth part of the objective function is a partly smooth mapping. The results are established using a strict complementary condition and a connection between our algorithm and the standard Riemannian trust-region method. Preliminary numerical results demonstrate that the approach performs promisingly on a class of $\ell_1$-optimization problems.

\appendix

\bibliographystyle{abbrv}
\bibliography{references}

\clearpage

\appendix

\section{Some supplementary examples}
\label{appendix 1}

\newtheorem{example-tau}[lem: partial psi]{Example}
\begin{example-tau} [Expressing $F_{\mathrm{nor}}^{\Lambda}(\tau(x))$] We consider the following examples:
\label{example-tau}
\begin{itemize}
\item[{(i)}] Group lasso: In this setting, $X=(X_1,\cdots,X_{n_2})$ is a matrix in $\mathbb{R}^{n_1\times n_2}$, $\psi=f+\varphi$, $\varphi(X)=\sum_{i=1}^{n_2}\left\Vert X_i\right\Vert_2$, and $\Lambda=\mathrm{diag}(\lambda_1 I_{n_1},\lambda_2 I_{n_1},\cdots, \lambda_{n_2} I_{n_1})$, where we understand $\Lambda$ by viewing $X$ as a vector $X=(X_1^T,\cdots,X_{n_2}^T)^T$. In this case, we obtain 
\begin{equation*}
    \partial\varphi(X)=\left\{W=(W_1,\cdots,W_{n_2})\in\mathbb{R}^{n_1\times n_2}: W_i\begin{cases}
    \in B_1(0) & \text{if } X_i=0,\\
    =X_i/\Vert X_i\Vert & \text{if } X_i\neq 0,
    \end{cases}\quad \forall~i=1,2,\cdots, n_2\right\},
\end{equation*}
and
\begin{equation*}
    F_{\mathrm{nor}}^{\Lambda}(\tau(X))_i=\begin{cases}
    \nabla f(X)_i+X_i/\left\Vert X_i\right\Vert & \text{if } X_i\neq 0,\\
    \nabla f(X)_i-\mathbf{P}_{B_1(0)}\left(\nabla f(X)_i\right), & \text{if }X_i=0,
    \end{cases} \quad \forall~i=1,2,\cdots,n_2.
\end{equation*}
\item[{(ii)}] $\ell_\infty$-optimization: $\psi=f+\varphi$, $\varphi(x)=\Vert x\Vert_{\infty}$ and $\Lambda=\lambda I$. Using the dual characterization $\|x\|_\infty = \max_{\|y\|_1 \leq 1}{x^Ty}$, we have
\begin{equation*}
\partial\varphi(x)= \{w\in\mathbb{R}^n:\Vert w\Vert_1\leq 1,\ w^T x = \|x\|_\infty \},
\end{equation*}
and $F_{\mathrm{nor}}^{\Lambda}(\tau(x))$ can be computed using \eqref{comp F tau}.
\end{itemize}
\end{example-tau}

\newtheorem{example-gamma}[lem: partial psi]{Example}
\begin{example-gamma}[Calculating $\Gamma(x)$]\label{example-gamma} We have:
\begin{itemize}
\item[{(i)}] Group lasso: For $X=(X_1,\cdots,X_{n_2})\in\mathbb{R}^{n_1\times n_2}$, $\psi=f+\varphi$, and $\varphi(X)=\sum_{i=1}^{n_2}\Vert X_i\Vert_2$, it holds that $\Gamma(X)=\min\{\Vert X_i\Vert: X_i\neq 0\}$.

\item[{(ii)}] $\ell_{\infty}$-optimization: Let us set $\psi=f+\varphi$, $\varphi(x)=\Vert x\Vert_{\infty}$, and $S:=\{i\in\{1,2,\cdots,n\}\mid |x_i|\neq\Vert x\Vert_{\infty} \}$. Then, it holds that
\begin{equation*}
    \Gamma(x)=\begin{cases}
    \Vert x\Vert_{\infty}-\max_{i\in S}|x_i| & S\neq \emptyset,\\
    2\Vert x\Vert_{\infty} & S=\emptyset.
    \end{cases}
\end{equation*}
\end{itemize}
\end{example-gamma}

\newtheorem{example-truncation}[lem: partial psi]{Example}
\begin{example-truncation}[Truncation]\label{example-truncation} We consider some additional functions that can be truncated:
\begin{itemize}
\item[{(i)}] Group lasso: For $X=(X_1,\cdots,X_{n_2})\in\mathbb{R}^{n_1\times n_2}$, $\psi=f+\varphi$, and $\varphi(X)=\sum_{i=1}^{n_2}\Vert X_i\Vert_2$, we can set $S_i=\{X\in\mathbb{R}^{n_1\times n_2}: \mathrm{card}\{j=1,2,\cdots,n_2\mid X_j=0\}\geq i\}$ for $i=1,2,\cdots,n_2$, $m=n_1$, $\delta=+\infty$, $\kappa=\sqrt{n_2}$, and $T(X,a)\in\mathbb{R}^{n_1\times n_2}$ is defined column-wise
\begin{equation*}
    T(X,a)_j={\mathbbm 1}_{\Vert \cdot \Vert \geq a}(X_j)\cdot X_j \quad j=1,2,\cdots,n_2.
\end{equation*}

\item[{(ii)}] $\ell_{\infty}$-optimization: For $\psi=f+\varphi$ and $\varphi(x)=\Vert x\Vert_{\infty}$, we can set $S_i=\{x\in\mathbb{R}^n\mid \mathrm{card}\{j=1,2,\cdots,n\mid x_j=\Vert x\Vert_{\infty}\}\geq i+1\}$ for $i=0,1,\cdots,n-1$, $S_n=\{0\}$, $m=n$, $\delta=+\infty$, and $\kappa=\sqrt{n}$. As for $T(x,a)\in\mathbb{R}^n$, if $x\in S_{n-1}$, it is defined via
\begin{equation*}
    T(x,a)={\mathbbm 1}_{\Vert \cdot \Vert_{\infty}\geq \frac{a}{2}}(x) \cdot x;
\end{equation*}
otherwise it is defined component by component via 
\begin{equation*}
    T(x,a)_j=x_j+{\mathbbm 1}_{|\cdot|> \Vert x\Vert_{\infty}- a}(x_j)\mathrm{sgn}(x_j)(\Vert x\Vert_{\infty}-|x_j|),\quad j=1,2,\cdots,n.
\end{equation*}
\end{itemize}
\end{example-truncation}

\section{The Solution of the Trust-region Subproblem}
\label{solution_sub}
In this section, we briefly discuss how to recover a solution $s^k$ satisfying \eqref{model reduce0} and \eqref{model reduce}. If the second-order information $B^k$ is symmetric, the subproblem \eqref{sub-problem} coincides with the classical trust-region subproblem and can be solved using standard methods, such as CG-Steihaug method \cite[Algorithm 7.2]{Num-Opti-Nocedal-06}. However, due to the nonsmoothness, the second-order information $B^k$ might be asymmetric. For example, the Jabobian of $g(x)=F_{\mathrm{nat}}^{\Lambda}(x)$ might be asymmetric. In this case, we can simply replace $B^k$ with its symmetrized version, $\frac{1}{2}[B^k+(B^k)^T]$, and then employ the CG-Steihaug method. 

If the matrix $B^k$ is positive semidefinite (but probably non-symmetric), i.e., $\langle h, B^k h\rangle\geq 0$ for all $h\in\mathbb{R}^n$, we can still solve \eqref{sub-problem} without symmetrization. We first choose a suitable regularization parameter $t_k\geq 0$ such that 

\begin{equation}
\frac{1}{2}h^T B^k h+t_k\Vert h\Vert^2\geq \lambda_1 \Vert h\Vert^2 \quad \forall\ h\in\mathbb{R}^n \quad \text{ and } \quad \Vert B^k+t_k I\Vert\leq \lambda_2,\label{lambda_sub}
\end{equation}
where $\lambda_1,\lambda_2$ are chosen positive constants which do not depend on $k$. We then consider the linear system
\begin{equation}
(B^k+t_k I)p=-g^k \label{Newton-sys0}
\end{equation}
and solve it to get an approximate solution $p^k$ satisfying
\begin{equation}
(B^k+t_k I)p^k=-g^k+r^k \quad \text{and} \quad \Vert r^k\Vert\leq\frac{\lambda_1}{2(\lambda_1+\lambda_2)}\Vert g^k\Vert,\label{appro-Newton-sys0}
\end{equation}
where $r^k$ is the residual. Finally, we project $p^k$ onto the trust region, i.e.,
\begin{equation}
s^k=\min\{\Delta_k,\Vert p^k\Vert\}{\bar p^k}.\label{proj_s}
\end{equation}
The next lemma proves that $s^k$ given by \eqref{Newton-sys0} and \eqref{proj_s} satisfies the condition \eqref{model reduce0} for some $\gamma_1,\gamma_2>0$.

\newtheorem{lem: model reduce}[lem: partial psi]{Lemma} 
\begin{lem: model reduce} Suppose that $B^k$ is positive semidefinite and that \eqref{lambda_sub} holds for some constants $\lambda_1,\lambda_2 >0$. Then condition \eqref{model reduce0} holds for $\gamma_1=\frac{\lambda_1}{2\lambda_2}$ and $\gamma_2=\frac{\lambda_1+2\lambda_2}{2\lambda_2(\lambda_1+\lambda_2)}$ if $s^k$ is given by \eqref{Newton-sys0} and \eqref{proj_s}.
\end{lem: model reduce}

\begin{proof}
Due to
\begin{equation*}
\begin{split}
\Vert p^k\Vert &= \Vert (B^k+t_k I)^{-1}(g^k-r^k)\Vert=\frac{\|g^k - r^k\|}{\Vert g^k-r^k\Vert/\left\Vert (B^k+t_k I)^{-1}(g^k-r^k)\right\Vert} \\ &\geq \frac{1}{\displaystyle\max_{h\in\mathbb{R}^n} \frac{\Vert (B^k+t_k I)h\Vert}{\Vert h\Vert}}\Vert g^k-r^k\Vert\geq \frac{1}{\lambda_2}\Vert g^k-r^k\Vert\geq \frac{\lambda_1+2\lambda_2}{2\lambda_2(\lambda_1+\lambda_2)}\Vert g^k\Vert
\end{split}
\end{equation*}
and utilizing the positive semidefiniteness of $B^k$, we have: 
\begin{align*}
m(0)-m(s^k)&\geq \frac{\min\{\Delta_k,\Vert p^k\Vert\}}{\Vert p^k\Vert}[m(0)-m(p^k)] \\
&=\frac{\min\{\Delta_k,\Vert p^k\Vert\}}{\Vert p^k\Vert}\left[ (p^k)^T(B^k+t_k I)p^k-(r^k)^T p^k-\frac{1}{2}(p^k)^T B^k p^k\right] \\
&\geq \frac{\min\{\Delta_k,\Vert p^k\Vert\}}{\Vert p^k\Vert}\left[ \frac{1}{2}(p^k)^T B^k p^k+t_k\Vert p^k\Vert^2-\frac{\lambda_1}{2(\lambda_1+\lambda_2)}\Vert g^k\Vert\Vert p^k\Vert\right] \\
&\geq \lambda_1\Vert p^k\Vert \min\{\Delta_k,\Vert p^k\Vert\}-\frac{\lambda_1}{2(\lambda_1+\lambda_2)}\Vert g^k\Vert\min\{\Delta_k,\Vert p^k\Vert\}\\
&\geq \frac{\lambda_1}{2\lambda_2} \Vert g^k\Vert\min\left\{\Delta_k,\frac{\lambda_1+2\lambda_2}{2\lambda_2(\lambda_1+\lambda_2)}\Vert g^k\Vert\right\}.
\end{align*}
Thus, \eqref{model reduce0} is satisfied for $\gamma_1=\frac{\lambda_1}{2\lambda_2}$ and $\gamma_2=\frac{\lambda_1+2\lambda_2}{2\lambda_2(\lambda_1+\lambda_2)}$.
\end{proof}

As for the other condition \eqref{model reduce}, we can simply set
\begin{equation}
s^k=\begin{cases}
    \text{the solution given by \eqref{Newton-sys0} and \eqref{proj_s}} & \text{if } \Delta_k\geq \zeta,\\
    s^k_C & \text{if } \Delta_k< \zeta
    \end{cases} \quad \text{and} \quad \ell(\Delta)=\begin{cases}0 & \text{if } \Delta<\zeta,\\ 1 & \text{if } \Delta\geq \zeta, \end{cases}
\end{equation}
where $\zeta > 0$ is a constant. We immediately obtain \eqref{model reduce}. 

\end{document}